\documentclass[oneside]{amsart}
\usepackage[a4paper, tmargin=1in, bmargin=1in]{geometry}
\usepackage[english]{babel}

\usepackage{amsmath}
\usepackage{graphicx}
\usepackage{graphics}
\usepackage{diagbox}
\usepackage{tikz-cd}

\usepackage{amsfonts}
\usepackage[mathcal]{euscript}
\usepackage{bbm}
\usepackage[colorlinks=true, allcolors=blue, linktocpage=true]{hyperref}
\usepackage{setspace}
\usepackage{amsthm}
\usepackage{cleveref}
\usepackage{float}
\usepackage{quiver}
\usepackage{mathtools}
\DeclareMathAlphabet{\mathpzc}{OT1}{pzc}{m}{it}
\usepackage{enumitem}
\setlist{itemsep=1pt}
\usepackage[style=alphabetic,natbib=true,firstinits=true,doi=false,isbn=false,url=false]{biblatex}
\newcommand{\curly}{\mathrel{\leadsto}}

\tikzcdset{every label/.append style = {font = \footnotesize},
every matrix/.append style={nodes={font=\small}}
}

\newcommand{\langl}{\begin{picture}(4.5,7)
\put(1.1,2.5){\rotatebox{60}{\line(1,0){5.5}}}
\put(1.1,2.5){\rotatebox{300}{\line(1,0){5.5}}}
\end{picture}}
\newcommand{\rangl}{\begin{picture}(4.5,7)
\put(.9,2.5){\rotatebox{120}{\line(1,0){5.5}}}
\put(.9,2.5){\rotatebox{240}{\line(1,0){5.5}}}
\end{picture}}

\theoremstyle{definition}
\newtheorem{thm}{Theorem}
\newtheorem{lemma}[thm]{Lemma}
\newtheorem{defn}[thm]{Definition}
\newtheorem{prop}[thm]{Proposition}

\newtheorem{conjecture}[thm]{Conjecture}
\newtheorem{remark}[thm]{Remark}
\numberwithin{thm}{subsection}

\newcommand{\sheaf}[0]{\mathcal{O}}
\newcommand{\sheaff}[0]{\mathcal{F}}
\newcommand{\db}[1]{\textbf{D}^{b}(#1)}
\newcommand{\dperf}[1]{\textbf{D}^{\operatorname{perf}}(#1)}
\newcommand{\sod}[1]{\langl #1 \rangl}
\newcommand{\gens}[1]{\langl #1 \rangl}
\newcommand{\pp}{\mathbb{P}^1}
\newcommand{\bl}[1]{\operatorname{Bl}_{\langle #1 \rangle}}
\newcommand{\bll}[1]{\operatorname{Bl}_{#1}}

\newcommand{\base}[1]{\pi_{B}^{*}\sheaf(#1)}

\newcommand{\mut}[1]{#1^{\prime}}
\newcommand{\mutt}[1]{#1^{\prime\prime}}
\newcommand{\fatpoint}{\mathbb{C}[s]/s^2}

\newcommand{\red}{\operatorname{red}}
\newcommand{\term}[1]{\textsf{#1}}
\newcommand{\cpinf}{c\mathbb{P}^{\infty}}
\newcommand{\pinf}[1]{\mathbb{P}^{\infty, #1}}

\newcommand{\kr}{\mathcal{K}r}
\newcommand{\exobj}{\mathcal{E}}
\newcommand{\restrict}[1]{\vert_{#1}}
\newcommand{\Hom}{\operatorname{Hom}}

\newcommand{\Ext}{\operatorname{Ext}}

\newcommand{\Extcomplex}{\operatorname{Ext}^{\bullet}}

\newcommand{\sh}[1]{\mathcal{#1}}
\newcommand{\pbundle}{\mathbb{P}}

\newcommand{\curve}{\mathfrak{C}}
\newcommand{\curvebase}{\curve_0}
\newcommand{\fpp}{\mathfrak{Y}}
\newcommand{\fppred}{\fpp_{\red}}

\newcommand\todo[1]{}

\begin{document}

\title{Categorical Absorption of a Non-Isolated Singularity}
\author{Aporva Varshney}
\begin{abstract}
   We study an example of a projective threefold with a non-isolated singularity and its derived category. 
   The singular locus can be locally described as a line of surface nodes compounded with a threefold node at the origin.
   We construct a semiorthogonal decomposition where one component absorbs the singularity in the sense of Kuznetsov--Shinder, and the other components are equivalent to the derived categories of smooth projective varieties.
   The absorbing category is seen to be closely related to the absorbing category constructed for nodal varieties by Kuznetsov--Shinder, reflecting the geometry of the singularity.
   We further show that the semiorthogonal decomposition is induced by one on a geometric resolution, and briefly consider the properties of the absorbing category under smoothing.
\end{abstract}
\maketitle

\tableofcontents

\section{Introduction}

Studying the bounded derived category of coherent sheaves $\db{X}$ of a variety $X$ over a field $k$ can often reveal information about the geometry of $X$.
One way in which we may study the structure of $\db{X}$ is by considering how it can be broken down into smaller pieces:
a \textit{semiorthogonal decomposition} 
\[
\db{X} = \sod{\sh{A}_1, \dots, \sh{A}_n}.
\]
is an ordered collection of full triangulated subcategories $\sh{A}_1, \dots, \sh{A}_n$ generating $\db{X}$, such that `there are no morphisms from right to left.'
Precisely, for $B \in \sh{A}_j$ and $A \in \sh{A}_i$ with $j > i$, $\Hom(B, A) = 0$.

Recently, there has been interest in studying $\db{X}$ when $X$ is a singular variety.
One approach we may take is to consider the full triangulated subcategory of perfect complexes, $\dperf{X}$, i.e. the complexes quasi-isomorphic to a bounded complex of locally free sheaves.
When $X$ is smooth, $\dperf{X}$ coincides with $\db{X}$, since in this setting all coherent sheaves have a finite resolution by locally free sheaves, and so we can think of $\dperf{X}$ as capturing smoothness.
It then makes sense to define the \textit{singularity category} \citep{orlovdsing}:
\[
\mathbf{D}_{\operatorname{sg}}(X) := \db{X} / \dperf{X}.
\]
This category exhibits locality in the sense that it depends only on an open neighbourhood of the singular locus.

On the other hand, one can take a more global approach to studying the singularity.
Suppose $X$ is a projective variety with a single isolated singularity. 
Then, one can consider whether the singularity is captured at the level of the derived category by asking if there exists a semiorthogonal decomposition
\[
\db{X} = \sod{ \mathcal{A}, \mathcal{B} }
\]
where $\mathcal{B} \subset \dperf{X}$, so that $\mathcal{A}$ is, in some sense, responsible for the singularity.
Following \citep{kuznetsov-shinder}, we say that $\sh{A}$ provides a \textit{categorical absorption of singularities}.
We can motivate this further as follows: there does \textit{not} exist a semiorthogonal decomposition of the form
\[
\db{X} = \sod{\mathbf{D}_{\operatorname{sg}}(X), \dperf{X}}
\]
but we can try to approximate it by asking that $\mathcal{B} \subset \dperf{X}$ is as large as possible.
Examples studied in the literature also have the property that $\sh{A}$ is equivalent to $\db{A}$ for some dg-algebra $A$ with finite dimensional cohomology.

The existence of a categorical absorption does not depend only on the singularity but also on the global structure of the variety; indeed, there are several obstructions as a result of this.
Nevertheless, we briefly review several instances where such a semiorthogonal decomposition has been found.
For the rest of the section we work over $k = \mathbb{C}$.

In \citep{kks}, the authors study surfaces with cyclic quotient singularities.
These include the singularities of the weighted projective space $\pbundle(1, 1, n)$, which are shown to have a semiorthogonal decomposition:
\[
\db{\pbundle(1, 1, n)} = \sod{\db{\mathbb{C}[z_1, \dots, z_{n - 1}]/(z_1, \dots, z_{n - 1})^2}, \db{\mathbb{C}}, \db{\mathbb{C}}}
\]
where $\deg z_i = 0$.
For surfaces with type $A_n$ singularities, it is shown that the singularity is often absorbed by the algebra $\mathbb{C}[w]/w^{n + 1}$ with $\deg w = 0$, with a simple example being the blow-up of a smooth surface $S$ in a non-reduced point $\operatorname{Spec} \mathbb{C}[w]/w^{n+1} \subset S$.
A semiorthogonal decomposition is produced for the threefold $\pbundle(1, 1, 1, 3)$ in \citep{kawamata}:
\[
\db{\pbundle(1, 1, 1, 3)} = \sod{\db{R}, \db{\mathbb{C}}, \db{\mathbb{C}}}
\]
where this time $R$ is an algebra with $\dim R = 45$.

The most relevant case for us will be categorical absorption of nodal varieties.
Nodal threefolds are considered in various incarnations in \citep{kawamata2}, \citep{pavic2021derived} and \citep{xie2023nodal}.
We follow the approach detailed in \citep{kuznetsov-shinder}.
To illustrate this theory, let $X$ be a nodal genus $0$ curve consisting of curves $E \cong Y \cong \pp$ intersecting transversely at a point.
The contraction of the branch $E$ through the map $f : X \to Y$ yields the semiorthogonal decomposition:
\[
\db{X} = \sod{\ker \textbf{R}f_{*}, \textbf{L}f^{*}\db{\pp}}
\]
where we have that $\textbf{L}f^*$ is fully faithful.
Absorption is then given by the category $\ker \textbf{R}f_{*} = \gens{\sheaf_{E}(-1)}$, and denoting $P := \sheaf_{E}(-1)$, one finds that:
\begin{align*}
    & \Extcomplex(P, P) \cong \mathbb{C}[\theta] \\
    & \ker \textbf{R}f_{*} \simeq \gens{P} \simeq \db{\mathbb{C}[r]/r^2}
\end{align*}
where $\deg \theta = 2$ and $\deg r = -1$.
We say that $P$ is a $\pinf{2}$-object, since $\Extcomplex(P, P)$ resembles the cohomology ring of $\mathbb{C}\mathbb{P}^{\infty}$.
It was shown in \textit{loc. cit.} that $\pinf{2}$-objects can be used to absorb singularities on nodal varieties of odd dimension.
In contrast, in even dimensions the absorbing category is instead given by $ \db{\mathbb{C}[w]/w^2} \simeq \gens{P}$ with $\deg w = 0$, and $\Extcomplex(P, P) = \mathbb{C}[\epsilon]$ with $\deg \epsilon = 1$. 
In this case, $P$ is called a $\pinf{1}$-object.

Perhaps the most interesting part of this story is a property of $\pinf{2}$-objects known as \textit{deformation absorption}: given a $\pinf{2}$ object on a variety, the pushforward into a smoothing of the variety generates an admissible subcategory of the smoothing. 
Essentially, this tells us that the $\pinf{2}$ object vanishes from the semiorthogonal decomposition when smoothing the variety.
In the example above, this is exemplified by considering a smoothing of $X$ to $\pp$.
This powerful property was used in \citep{ks-fano} to study relationships between the Kuznetsov components of smooth Fano threefolds, by first degenerating a threefold to a nodal variety and then taking a resolution.

\subsection{Results and structure}

The cases in the literature have so far focused on isolated singularities. 
Presently, we turn our attention to a complex threefold with a non-isolated compound $A_2$ singularity, constructed as follows.
We first consider the curve
\[
\curve := \{x = s^2u = 0\} \subset \pp_{x:y} \times \pp_{s:t} \times \pp_{u:v}
\]
and then perform a blow-up:
\[
X := \bll{\curve}(\pp)^3.
\]
Then, $X$ is locally given by the equation $\{xb = s^2u\} \subset \mathbb{A}^4_{x, s, u, b}$, so we can think of the singularity as a family of surface nodes compounded with a threefold node at the origin.
In this article, we essentially perform a series of calculations checking that there is a categorical absorption which is as nicely behaved as one could expect.

\begin{thm}[=\ref{thm:curveabsorption}]
    The curve $\curve$ and threefold $X$ admit a categorical absorption, where the absorbing category is generated by an object $P$ with the property
    \[
    \Extcomplex(P, P) \cong \mathbb{C}[\epsilon, \theta]
    \]
    where $\deg \epsilon = 1$ and $\deg \theta = 2$.
\end{thm}

We call an object with this self-Ext algebra a \term{compound} $\mathbb{P}^{\infty}$, or $\cpinf$, object, to reflect the idea that it is simultaneously a $\pinf{1}$ and $\pinf{2}$ object.

The work in \citep{ks-fano} shows us that it can be useful to consider the structure of a resolution and of a smoothing; consequently, in \Cref{section:catres} we show that the semiorthogonal decomposition produced for $X$ is induced by one on a specific choice of resolution of singularities $\widetilde{X}$.
In particular, our resolution has as exceptional divisor a surface and a curve intersecting transversely.
The key geometric input which allows us to perform computations in the derived category of $\widetilde{X}$ is the following proposition.
\begin{prop}[=\ref{prop:altres}]
    The space $\widetilde{X}$ can be obtained from $(\pp)^3$ as a series of blow-ups along smooth rational curves.
\end{prop}

Using this, we obtain a full exceptional collection on $\widetilde{X}$.
However, this collection is not compatible with the contraction $\sigma : \widetilde{X} \to X$, in the sense that the pushforward under $\sigma$ does not recover the absorption constructed on $X$. 
To overcome this, we perform a series of mutations and explicitly find an admissible category $\widetilde{\mathcal{D}} \subset \db{\widetilde{X}}$ that serves as a `categorical resolution' for our absorbing category.
\begin{thm}[=\ref{thm:resolution}]
    Let $\sigma : \widetilde{X} \to X$ be the resolution of singularities map.
    There is an admissible subcategory $\widetilde{\mathcal{D}} \subset \db{\widetilde{X}}$ such that the absorbing category $\mathcal{D} \subset \db{X}$ is equivalent to the Verdier quotient
    \[
    \mathcal{D} \simeq \widetilde{\mathcal{D}} / (\ker \mathbf{R}\sigma_{*} \cap \widetilde{\sh{D}}).
    \]
\end{thm}

The category $\widetilde{\mathcal{D}}$ is generated by an exceptional collection of length $4$, and under the pushforward $\mathbf{R}\sigma_{*}$ all the generating objects are mapped to the $\cpinf$ object $P$. 
We investigate the algebra structure of this collection.
\begin{thm}[=\ref{thm:resolutionmorphisms}]
    The morphisms in the exceptional sequence
    \[
    \widetilde{\mathcal{D}} = \sod{\sh{E}_0^{\prime\prime}, \sh{F}_0^{\prime\prime}, \sh{E}_1, \sh{F}_1}
    \]
    are described as follows:
    \begin{enumerate}
        \item $\Extcomplex(\sh{E}_0^{\prime\prime}, \sh{F}_0^{\prime\prime}) \cong \Extcomplex(\sh{E}_1, \sh{F}_1) \cong \mathbb{C} \oplus \mathbb{C}[-1]$
        \item $\Extcomplex(\sh{F}_0^{\prime\prime}, \sh{E}_1) \cong \mathbb{C} \oplus \mathbb{C}[-1]$
        \item $\Extcomplex(\sh{F}_0^{\prime\prime}, \sh{F}_1) \cong \Extcomplex(\sh{E}_0^{\prime\prime}, \sh{E}_1) \cong \mathbb{C} \oplus \mathbb{C}[-1] \oplus \mathbb{C}^{\oplus 2}[-2]$
        \item $\Extcomplex(\sh{E}_0^{\prime\prime}, \sh{F}_1) \cong \mathbb{C} \oplus \mathbb{C}[-1] \oplus \mathbb{C}^{\oplus 2}[-2] \oplus \mathbb{C}^{\oplus 2}[-3]$
    \end{enumerate}
\end{thm}

To compute how morphisms compose, we note that there is an isomorphism of graded vector spaces between the Ext-spaces in the theorem above and the truncated self-Ext-spaces of the $\cpinf$ object.
\begin{align*}
    & \Extcomplex(\sh{E}_0^{\prime\prime}, \sh{F}_0^{\prime\prime}) \cong \Extcomplex(\sh{E}_1, \sh{F}_1) \cong \Extcomplex(\sh{F}_0^{\prime\prime}, \sh{E}_1) \cong \Ext^{\leq 1}(P, P) \\
    & \Extcomplex(\sh{E}_0^{\prime\prime}, \sh{E}_1)  \cong \Extcomplex(\sh{F}_0^{\prime\prime}, \sh{F}_1) \cong \Ext^{\leq 2}(P, P) \\
    & \Extcomplex(\sh{E}_0^{\prime\prime}, \sh{F}_1) \cong \Ext^{\leq 3}(P, P)
\end{align*}
This indicates that the morphisms between objects of $\widetilde{\mathcal{D}}$ should be seen as lifting some map $P \to P[i]$.
We prove the following.
\begin{prop}[=\ref{prop:algebrastructure}]
    All compositions in $\widetilde{\mathcal{D}}$ are non-zero.
\end{prop}
In particular, omiting compositions in non-zero degrees, we can visualise the categorical resolution as the following quiver.
\[\begin{tikzcd}[row sep=large,column sep=large]
	{\mathcal{E}_0^{\prime\prime}} && {\mathcal{F}_0^{\prime\prime}} \\
	\\
	{\mathcal{E}_1} && {\mathcal{F}_1}
	\arrow["{t_0}"{description}, curve={height=-6pt}, from=1-1, to=1-3]
	\arrow["{\epsilon_0}"{description}, curve={height=6pt}, dashed, from=1-1, to=1-3]
	\arrow["{s_0}"{description}, curve={height=-6pt}, from=1-1, to=3-1]
	\arrow["{t_1}"{description}, curve={height=6pt}, from=1-3, to=3-1]
	\arrow["{\epsilon_1}"{description}, curve={height=-6pt}, dashed, from=1-3, to=3-1]
	\arrow["{\theta_1}"{description}, curve={height=-6pt}, dotted, from=1-3, to=3-3]
	\arrow["{t_2}"{description}, curve={height=-6pt}, from=3-1, to=3-3]
	\arrow["{\epsilon_2}"{description}, curve={height=6pt}, dashed, from=3-1, to=3-3]
	\arrow["{s_1}"{description}, curve={height=6pt}, from=1-3, to=3-3]
	\arrow["{\theta_0}"{description}, curve={height=6pt}, dotted, from=1-1, to=3-1]
\end{tikzcd}\]
Here, solid arrows indicate degree $0$ maps, dashed arrows indicate degree $1$ maps and dotted arrows indicate degree $2$ maps.

Finally, we briefly consider the properties of our $\cpinf$ objects with respect to smoothings of the curve $\curve$ in \Cref{section:smoothing}.
\begin{prop}[=\ref{prop:defabs}]
    Let $\mathfrak{X}$ be the smoothing of $\curve$ obtained by blowing up $\pp \times \pp$ in a fattened point $\operatorname{Spec} \mathbb{C}[s]/s^2$.
    Then the pushforward of the $\cpinf$ object $P \in \db{\curve}$ into $\mathfrak{X}$ generates an admissible subcategory of $\db{\mathfrak{X}}$, and hence provides a deformation absorption with respect to the smoothing $\mathfrak{X}$.
\end{prop}

\subsubsection*{Notation and Conventions}
\begin{enumerate}
    \item We work over the base field $\mathbb{C}$.
    \item For a dg algebra $A$ with finite dimensional cohomology, we denote by $\mathbf{D}(A)$ the unbounded derived category of dg-modules over $A$, and $\db{A} \subset \mathbf{D}(A)$ is the subcategory of dg-modules with finite dimensional total cohomology.
    \item Tensor products and pushforward and pullback functors between derived categories are implicitly derived.
    \item Given a scheme $X$ with closed subscheme $Z \subset X$ and a birational map $Y \to X$, we denote by $\overline{Z}$ the strict transform of $Z$ in $Y$.
\end{enumerate}

\subsubsection*{Acknowledgements}
The author would like to thank his PhD supervisor Ed Segal for invaluable guidance throughout the project, and the anonymous reviewer for helpful comments which improved the paper.

\section{Preliminaries}

\subsection{Derived categories}
In this section we recap general results in the theory of derived categories which we will use freely.

\subsubsection{Mutations of semiorthogonal decompositions}\label{sec:mutation}
Given a semiorthogonal decomposition of a triangulated category where each component is admissible
\[
\sh{T} = \sod{\sh{A}_1, \dots, \sh{A}_n}
\]
it is possible to reorder the components to obtain a new semiorthogonal decomposition.
The first method is to perform a \textit{left mutation} of $\sh{A}_{i + 1}$ over $\sh{A}_{i}$:
\[
\sh{T} = \sod{\sh{A}_1, \dots, \sh{A}_{i - 1}, \mathbb{L}_{\sh{A}_i}\sh{A}_{i + 1}, \sh{A}_i, \sh{A}_{i + 2}, \dots, \sh{A}_n}
\]
and the second is to perform a \textit{right mutation} of $\sh{A}_{i}$ over $\sh{A}_{i + 1}$:
\[
\sh{T} = \sod{\sh{A}_1, \dots, \sh{A}_{i - 1}, \sh{A}_{i + 1}, \mathbb{R}_{\sh{A}_{i + 1}}\sh{A}_{i}, \sh{A}_{i + 2}, \dots, \sh{A}_n}.
\]
In each case, the mutation is an equivalence of categories:
\begin{align*}
    & \mathbb{L}_{\sh{A}_i}\sh{A}_{i + 1} \simeq \sh{A}_{i + 1} \\
    & \mathbb{R}_{\sh{A}_{i + 1}}\sh{A}_{i} \simeq \sh{A}_{i}
\end{align*}
We omit the construction, but note that when $\sh{A}_i$ and $\sh{A}_{i + 1}$ are generated by exceptional objects $\sh{E}_i$ and $\sh{E}_{i + 1}$ respectively, we have formulae for the left and right mutations through the following exact triangles \citep{b90}:
\begin{align}
    & \operatorname{RHom}(\sh{E}_i, \sh{E}_{i + 1}) \otimes \sh{E}_i \to \sh{E}_{i + 1} \to \mathbb{L}_{\sh{E}_i}
    (\sh{E}_{i + 1}) \label{eq:leftmutation} \\
    & \mathbb{R}_{\sh{E}_{i + 1}}(\sh{E}_i) \to \sh{E}_i \to \operatorname{RHom}(\sh{E}_i, \sh{E}_{i + 1})^{*} \otimes \sh{E}_{i + 1}. \label{eq:rightmutation}
\end{align}

We also note an easy consequence of Serre duality.
\begin{lemma}\label{serreduality}
If $X$ is Gorenstein with a semiorthogonal decomposition
\[
\db{X} = \sod{\sh{A}, \sh{B}}
\]
we have that
\[
\db{X} = \sod{\sh{B} \otimes \omega_{X}, \sh{A}}
\]
is also a semiorthogonal decomposition, where $\omega_X$ is the canonical bundle.
\end{lemma}

\subsubsection{Derived nine-lemma}
When working in the derived category, the ordinary nine-lemma in homological algebra no longer holds. 
Instead, we have the following \textit{derived} nine-lemma (see \citep[Lemma 2.6]{may} for a proof).

\begin{lemma}\label{lemma:derivednine}
    Suppose that the top left square commutes in the following diagram and the top two rows and left two columns are exact triangles.
    \[\begin{tikzcd}
	X & Y & Z & {X[1]} \\
	{X'} & {Y'} & {Z'} & {X'[1]} \\
	{X''} & {Y''} & {Z''} & {X''[1]} \\
	{X[1]} & {Y[1]} & {Z[1]} & {X[2]}
	\arrow[from=1-1, to=1-2]
	\arrow[from=1-2, to=1-3]
	\arrow[from=1-2, to=2-2]
	\arrow[from=2-2, to=3-2]
	\arrow[from=1-1, to=2-1]
	\arrow[from=2-1, to=2-2]
	\arrow[from=2-2, to=2-3]
	\arrow[from=2-1, to=3-1]
	\arrow["f"', dashed, from=3-1, to=3-2]
	\arrow["g"', dashed, from=3-2, to=3-3]
	\arrow["k", dashed, from=1-3, to=2-3]
	\arrow["{k'}", dashed, from=2-3, to=3-3]
	\arrow[from=1-3, to=1-4]
	\arrow[from=1-4, to=2-4]
	\arrow[from=2-4, to=3-4]
	\arrow["{k''}", dashed, from=3-3, to=4-3]
	\arrow[from=3-2, to=4-2]
	\arrow[from=3-1, to=4-1]
	\arrow[from=4-1, to=4-2]
	\arrow[from=4-2, to=4-3]
	\arrow[from=2-3, to=2-4]
	\arrow["h"', dashed, from=3-3, to=3-4]
	\arrow[from=3-4, to=4-4]
	\arrow[from=4-3, to=4-4]
\end{tikzcd}\]    
    Then there exists an object $Z^{\prime\prime}$ and maps $f, g, h, k, k^{\prime}, k^{\prime\prime}$ such that the bottom right square anti-commutes, all other squares commute and all rows and columns form exact triangles.
\end{lemma}

\subsubsection{Derived category of a blow-up}

Suppose $X$ is a scheme and $Z \subset X$ is a locally complete intersection of codimension $c \geq 2$.
Consider the blow-up $\bll{Z}{X}$.
Then there is a semiorthogonal decomposition \citep{orlov93}:
\begin{equation}\label{eq:sodblowup}
    \db{\bll{Z}{X}} = \sod{ \Phi_{c - 1}(\db{Z}), \dots, \Phi_1(\db{Z}), \db{X}}
\end{equation}
where $\db{X}$ is embedded fully faithfully by pulling back, and the  functors $\Phi_i$ are fully faithful and constructed as follows.
The exceptional divisor is a projective bundle over $\pi : E \to Z$, so letting $j : E \to \bll{Z}{X}$ be the inclusion, we define
\begin{equation}\label{eq:blowupembedding}
    \Phi_{i}(\sh{F}) = j_{*}(\pi^{*}\sh{F} \otimes \sheaf(-i)).
\end{equation}

\subsubsection{The cone on the pull-push counit}

Let $i : E \to X$ be the inclusion of a Cartier divisor in a scheme $X$.
A useful computation tool will be the following description of the cone on the canonical adjunction morphism \citep{bo95}:
\begin{equation}\label{eq:counittriangle}
    i^{*}i_{*}\sh{F} \to \sh{F} \to \sh{F} \otimes \sh{N}_{E / X}^{\vee}[2].
\end{equation}

\subsubsection{Bondal--Orlov localization}\label{sec:boloc}
The Bondal--Orlov localization conjecture relates the derived categories of a singular variety $X$ and a resolution $\pi : \widetilde{X} \to X$.

\begin{conjecture} \citep[Conjecture 1.9]{efi20} \citep[Section 5]{bo02}
    If $\pi : \widetilde{X} \to X$ is a resolution of rational singularities, then $\pi_{*}$ induces an equivalence between $\db{X}$ and a Verdier quotient of $\db{\widetilde{X}}$:
    \[
    \widetilde{\pi}_{*} : \db{\widetilde{X}} / \ker \pi_{*} \to \db{X}.
    \]
\end{conjecture}

In fact the conjecture is true in several settings (see \citep{efi20, ps21}).
The most relevant situation for us will be when the resolution has one-dimensional fibres, in which case the conjecture was proven in \citep{bks}. 
An immediate consequence of this is that we may descend an appropriate semiorthogonal decomposition of $\db{\widetilde{X}}$ to one for $\db{X}$ \citep[Proposition 4.1]{kuznetsov-shinder}.
If
\[
\db{\widetilde{X}} = \sod{\widetilde{\sh{A}}, \widetilde{\sh{B}}}
\]
is a semiorthogonal decomposition such that $\ker \pi_{*} \subset \widetilde{\sh{A}}$, then we obtain a semiorthogonal decomposition
\[
\db{X} = \sod{\sh{A}, \sh{B}}
\]
with
\begin{align*}
    & \sh{A} \simeq \widetilde{\sh{A}} / (\ker\pi_{*} \cap \widetilde{\sh{A}}) \simeq \widetilde{\sh{A}} / \ker\pi_{*} \\
    & \sh{B} \simeq \widetilde{\sh{B}} / (\ker \pi_{*} \cap \widetilde{\sh{B}}) \simeq \widetilde{\sh{B}}
\end{align*}
where $\sh{A}, \sh{B}$ must generate $\db{X}$ because $\pi_{*}$ is essentially surjective.

\subsubsection{Spherical Objects}\label{sec:spherical}

An object $K \in \db{X}$, where $X$ is a complex smooth and projective variety, is called \textit{spherical} if \citep{huybrechts}:
\begin{enumerate}
    \item $K \otimes \omega_{X} \cong K$, and
    \item $\Hom(K, K[i]) \cong H^{i}(S^{\dim X}, \mathbb{C})$
\end{enumerate}
where $S^{\dim X}$ is the sphere of dimension equal to $\dim X$.
A spherical object induces an autoequivalence of the derived category by \textit{twisting} around the object.
For any object $\sh{E}$, the twist $T_{K}(\sh{E})$ is defined through the exact triangle
\[
\left( \bigoplus_{i} \Hom(K, \sh{E}[i]) \otimes K[-i] \right) \to \sh{E} \to T_{K}(\sh{E}).
\]
Being an autoequivalence, it has an inverse known as the dual twist, which is similarly described by an exact triangle.

\subsection{Projective bundles and ruled surfaces}\label{sec:ruledsurfaces}

We review some theory which will be essential for computations later on.
Fix a smooth variety $B$, and let $\sheaff$ and $\sh{L}$ be invertible sheaves.
Let $\sh{E} = \sh{L} \oplus \sheaff$.
Then there are two possible conventions for the corresponding projective bundle; in our case, we take the definition to be 
\[\pi : \mathbb{P}(\sh{E}) := \underline{\operatorname{Proj}}(\operatorname{Sym}(\sh{E})) \to B\] so that there is a map $\pi^{*}\sh{E} \to \sheaf(1)$.
We will write $\pbundle(\sh{E}) = \pbundle(\sh{L}_{x_0} \oplus \sh{F}_{x_1})$ to denote that the coordinates on the fibres are $[x_0 : x_1]$.
There is a section $S$ corresponding to the surjection $\sh{E} \to \sh{L}$; equivalently, it is given by the vanishing of the coordinate $x_1$.
The skyscraper sheaf of the section is defined by taking the cokernel of the composition $\pi^{*}\sheaff \to \pi^{*}\sh{E} \to \sheaf(1)$:
\[
\pi^{*}\sheaff \to \sheaf(1) \to \sheaf_{S}(1) \to 0
\]
and twisting by $\sheaf(-1)$
\[
\pi^{*}\sheaff \otimes \sheaf(-1) \to \sheaf \to \sheaf_{S} \to 0.
\]
Then the above sequence is the Koszul resolution of $\sheaf_S$, so it follows from the fact that $\sheaf(1)\restrict{S} \cong \sh{L}$ that the normal bundle of $S$ is given by $\sh{F}^{\vee} \otimes \sh{L}$.
Note that the first map is ``multiplication by $x_1$,'' and hence $x_1$ is a section of the line bundle $\pi^*\sheaff^{\vee} \otimes \sheaf(1)$.

We collect some results on ruled surfaces, which will arise frequently when we blow-up a curve $C \cong \pp$ in a smooth threefold $X$.
In this situation, we recall that the exceptional divisor is isomorphic to $\pbundle_{C}(\sh{N}_{C / X}^{\vee})$ and has normal bundle $\sheaf(-1)$.
Consider 
\[
\pi : \Sigma = \pbundle_{\pp}(\sheaf_x \oplus \sheaf(a)_y) \to \pp.
\]
Then, the divisor corresponding to $\sheaf(1)$ is given by the class of the section $H := V(x)$, and we denote
\begin{equation}\label{eq:ruledpicard}
\sheaf(a, b) := \sheaf(1)^{\otimes a} \otimes \pi^{*}\sheaf(b)    
\end{equation}
or, equivalently, it is the line bundle corresponding to the divisor $aH + bF$ where $F$ is a fibre.
Then the canonical bundle is given by the formula
\begin{equation}\label{eq:ruledcanonical}
    \omega_{\Sigma} \cong \sheaf(-2, -2 + a).
\end{equation}
If $\exobj = \sheaf \oplus \sheaf(a)$, then pushing forward the line bundle $\sheaf(d)$ to the base $\pp$ gives us:
\begin{equation}\label{eq:linebundlepushforward}
    \pi_{*}\sheaf(d) = 
\begin{cases}
    \operatorname{Sym}^d(\exobj) &\text{if $d \geq 0$}\\
    \left(\sheaf(a) \otimes \operatorname{Sym}^{-2 - d}(\exobj)\right)^{\vee}[-1] &\text{if $d < 0$}
\end{cases}
\end{equation}

\subsection{Categorical absorption for nodal varieties}

This section gives an overview of the theory of categorical absorption for nodal varieties, with \citep{kuznetsov-shinder} as our primary reference.
We will only concern ourselves with dimensions up to $3$, but note that the theory can be extended to arbitrary dimension.

\subsubsection{Categorical ordinary double points}\label{sec:odp}

We recall the definitions given earlier.
\begin{defn}
    An object $P \in \db{X}$ for a variety $X$ is called a $\pinf{q}$ object if
    \begin{align*}
    \Extcomplex(P, P) \cong \mathbb{C}[\theta] && \deg \theta = q
    \end{align*}
\end{defn}
\begin{defn}
    The \term{categorical ordinary double point of degree $p$} is $\db{\mathbb{C}[w]/w^2}$ with $\deg w = p$.
\end{defn}

The theory is split up by parity of $\dim X$: in all examples, $q = 1$ when $X$ is even dimensional, and $q = 2$ when $X$ is odd dimensional.

For the odd dimensional case the key example is the nodal curve considered earlier, which we can see as the vanishing locus \[X = V(su) \subset \pp_{s:t} \times \pp_{u:v}.\]
In this case, the pushforward of a line bundle $\sheaf(a)$ on either of the two branches gives a $\pinf{2}$ object. 
We have seen that we can consider contracting one of the branches $f : X \to \pp_{s:t}$; it is then a direct consequence of the fact that $f_{*}\sheaf_X \cong \sheaf_{\pp}$ that the pullback $f^{*}$ is fully faithful.
One can show that the null category $\ker f_{*}$ is admissible in this case, and that there is a semiorthogonal decomposition \citep[Proposition 6.15]{kuznetsov-shinder}
\[
\db{X} = \sod{\db{\mathbb{C}[w]/w^2}, \db{\pp}}
\]
with $\deg w = -1$.
Here, the absorbing category $\db{\mathbb{C}[w]/w^2}$ is generated by the pushforward of the $\sheaf(-1)$ sheaf on the branch being contracted.

The simplest (but still useful) example of a nodal even dimensional scheme is the fat point $\operatorname{Spec} \mathbb{C}[w]/w^2$.
In this case, the object $\mathbf{\mathbb{C}} \cong \sheaf / w$ is a $\pinf{1}$ object which generates $\db{\mathbb{C}[w]/w^2}$, where $\deg w = 0$.
More interestingly, we may consider the nodal surface $X := \mathbb{P}(1, 1, 2)$ \citep[Example 3.17]{kks}.
To construct a categorical absorption here, we consider the resolution of singularities given by
\[
\sigma : \mathbb{F}_2 \to \mathbb{P}(1, 1, 2)
\]
where $\pi : \mathbb{F}_2 \cong \pbundle(\sheaf \oplus \sheaf(-2)) \to \pp$ is the second Hirzebruch surface.
In this case the exceptional divisor $E$ is given by the section of the bundle which has self-intersection $(-2)$, and it can be shown that the kernel of $\sigma_{*}$ is generated by $\sheaf_E(-1)$.
Let $F$ be the class of a fibre, so that $\sheaf(F) \cong \pi^{*}\sheaf(1)$.
Then, $E$ is the divisor corresponding to $\sheaf(1)$ and so there is a full exceptional sequence:
\begin{align*}
\db{\mathbb{F}_2} & = \sod{\pi^{*}\db{\pp} \otimes \sheaf(-1), \pi^{*}\db{\pp}} \\ & =\sod{\sheaf(-E - 2F), \sheaf(-E - F), \sheaf(-F), \sheaf}.    
\end{align*}
One sees that there is a short exact sequence:
\[
0 \to \sheaf(-E - F) \to \sheaf(-F) \to \sheaf_E(-1) \to 0
\]
which implies that $\ker \sigma_{*} \subset \gens{\sheaf(-E - F), \sheaf(-F)}$.
An important feature here is that $\sheaf_E(-1)$ is a $2$-spherical object, and $\sheaf(-E - F)$ is the dual spherical twist of $\sheaf(-F)$ around $\sheaf_E(-1)$.
We note that the pushforwards $\sigma_{*}\sheaf$ and $\sigma_{*}\sheaf(-E-2F)$ remain exceptional in $\db{\mathbb{P}(1, 1, 2)}$, and hence $\sigma_{*}\gens{\sheaf(-E, -F), \sheaf(-F)}$ is a candidate for our categorical absorption.
We can compute, using \eqref{eq:linebundlepushforward}, that \[\Extcomplex(\sheaf(-E-F), \sheaf(-F)) \cong \mathbb{C} \oplus \mathbb{C}[-1].\]
Since the cone on the morphism $\sheaf(-E-F) \to \sheaf(-F)$ is in $\ker \sigma_{*}$, we find that \[\sigma_{*}\sheaf(-E-F) \cong \sigma_{*}\sheaf(-F).\]
The idea is that this object should be a $\pinf{1}$ object; intuitively, the degree zero morphism becomes an isomorphism and the degree $1$ morphism then gives rise to the self extension of the $\pinf{1}$ object.

These observations are formalised in the following theorems.
\begin{defn}\citep[\S 3.1]{kuznetsov-shinder}\label{def:gradedkronecker}
    The \term{graded Kronecker quiver of degree $q$} is the path dg-algebra of the following dg-quiver:
    \[\begin{tikzcd}
	\bullet && \bullet
	\arrow["0", curve={height=-6pt}, from=1-1, to=1-3]
	\arrow["q"', curve={height=6pt}, dashed, from=1-1, to=1-3]
\end{tikzcd}\]
    where there are two arrows of degrees $0$ and $q$ respectively, and the differentials are zero.
    The \term{graded Kronecker quiver category of degree $q$}, denoted $\kr_q$, is the derived category of this algebra.
\end{defn}

\begin{thm}\citep[Corollary 3.3]{kuznetsov-shinder}
    Let $\gens{\sh{E}, \sh{E}^{\prime}} \subset \db{\widetilde{X}}$ be an exceptional sequence, where $\widetilde{X}$ is a variety, such that $\Extcomplex(\sh{E}, \sh{E}^{\prime}) \cong k \oplus k[-q]$.
    Then $\gens{\sh{E}, \sh{E}^{\prime}} \simeq \kr_q$.
\end{thm}

\begin{thm} \citep[Proposition 3.7]{kuznetsov-shinder}
    Let $\gens{\sh{E}, \sh{E}^{\prime}} \simeq \kr_q$, and define the object $K$ using the distinguished triangle:
    \[
    K \to \sh{E} \to \sh{E}^{\prime}.
    \]
    Then, $K$ is a $1 + q$-spherical object, and there is an equivalence
    \begin{align*}
        \kr_q / \gens{K} \simeq \db{\mathbb{C}[w]/w^2} && \deg w = 1 - q
    \end{align*}
    and the image of $\sh{E}$ under $\sigma_{*} : \kr_q \to \kr_q / \gens{K}$ is a $\pinf{q}$ object.
\end{thm}

In the more general setting, the theorem can be applied as follows.
Suppose we have a crepant resolution $\sigma : \widetilde{X} \to X$ for a nodal surface or threefold $X$, and assume that $H^{>0}(\widetilde{X}, \sheaf) = 0$, so that line bundles are exceptional objects.
Then labeling the exceptional curve as $C$ we have a $d$-spherical object $K := \sheaf_C(-1)$, where $d = \dim X$.
One can show that $K$ generates $\ker \sigma_{*}$.
Now suppose $\sh{E}$ is an exceptional object such that $\dim \Extcomplex(\sh{E}, K) = 1$.
A simple way to find such an object is to look for a divisor which intersects the curve transversely at a point, and then take the associated line bundle.
Now, we can simply take the spherical twist of $\sh{E}$ around $K$; if $\Extcomplex(K, \sh{E})$ is concentrated in degree $n$, this takes the form:
\[
K[-n] \to \sh{E} \to \sh{E}^{\prime}
\]
Since a spherical twist is an autoequivalence, we see that $\sh{E}^{\prime}$ is an exceptional object, so we have an admissible subcategory $\gens{\sh{E}, \sh{E}^{\prime}}$ of $\db{\widetilde{X}}$ which contains $\ker \sigma_{*} = \gens{K}$.
It is then not hard to show, simply by taking long exact sequences for $\Ext$, that $\gens{\sh{E}, \sh{E}^{\prime}} \simeq \kr_{q}$ for $q = d - 1$, and hence taking the Verdier quotient gives us a $\mathbb{P}^{\infty, q}$ object generating a categorical ordinary double point of degree $1 - q$.
As a result, we are able to construct a categorical absorption for the nodal variety.

\subsubsection{Deformation absorption}\label{subsec:defabs}
Finally we discuss smoothings of nodal varieties.
Suppose that $f : \mathfrak{X} \to B$ is a flat projective morphism from a smooth space $\mathfrak{X}$ to a smooth curve $B$, with the central fibre isomorphic to $X$ and all other fibres smooth, so that $\mathfrak{X}$ gives a smoothing of $X$.
Now suppose that the nodal variety $X$ has an absorption by a $\pinf{2}$ object $P$:
\[
\db{X} = \sod{P, \prescript{\perp}{}{P}}.
\]
Then for any smoothing $i : X \xhookrightarrow{} \mathfrak{X}$, the object $i_{*}P$ can be shown to be exceptional \citep[Theorem 1.8]{kuznetsov-shinder}. 
In particular, we have a semiorthogonal decomposition
\[
\db{\mathfrak{X}} = \sod{i_{*}P, \sh{D}}
\]
which is said to be $B$-linear in the sense that there is a well-defined notion of base changing the semiorthogonal decomposition to a fibre of the map $f : \mathfrak{X} \to B$ \citep{basechange}.
Performing this base-change to a non-central fibre $\mathfrak{X}_{b}$, we find that the $i_{*}P$ component vanishes since it is only supported on the central fibre:
\[
\db{\mathfrak{X}_{b}} \simeq \sh{D}_{b}.
\]
In effect, we have a family of smooth and proper categories $\sh{D}_b$, and for $b = o$ we acquire a $\pinf{2}$ object to account for the singular behaviour.
We say that $\pinf{2}$ objects provide a \term{deformation absorption}, and furthermore, since this occurs for every smoothing, we say that the deformation absorption is \term{universal}.

\section{Compound \texorpdfstring{$\mathbb{P}^{\infty}$}{P-infinity} Objects} \label{section:cpinfobjects}

We begin by generalising the example of a nodal curve by trivially thickening one of the branches.
Consider 
\begin{equation}\label{eq:curve}
    \curve = \{s^2u = 0\} \subset \pp_{s:t} \times \pp_{u:v}
\end{equation}
which is a copy of $\curvebase := V(u) \cong \pp_{s:t}$ intersecting the thickened branch \[\fpp := V(s^2) \cong \pp_{u:v} \times \mathbb{C}[s]/s^2\] at a fat point.
There is a contraction $f : \curve \to \curvebase$, and it follows once again from the fact that $f_{*}\sheaf_{\curve} = \sheaf_{\curvebase}$ that there is a semiorthogonal decomposition
\[
\db{\curve} = \sod{\ker f_{*}, f^{*}\db{\curvebase}}.
\]
with $f^{*}$ fully faithful.
Let $\fppred = V(s)$ be the reduction of the fattened branch $\fpp$.
Denote by $\sh{P}(a)$ the pushforward of the line bundle $\sheaf(a)$ on $\fppred$ to $\curve$.
\begin{lemma}
    The kernel of $f_{*}$ is generated by $\sh{P}(-1)$.
\end{lemma}
\begin{proof}
    Firstly consider the projections $q : \fpp \to \mathbb{C}[s]/s^2$ and $p : \fpp \to \pp$.
    Then we have that
    \[
    \ker q_{*} = \gens{q^{*}\mathbf{k} \otimes p^{*}\sheaf(-1)}
    \]
    where we denote by $\mathbf{k} = (\mathbb{C}[s]/s^2)/s$ the generator of $\db{\mathbb{C}[s]/s^2}$.
    The pushforward of this sheaf to $\curve$ is precisely $\sh{P}(-1)$.

    To reduce to this case, we show that any object in the kernel of $f_*$ is a pushforward of an object in $\db\fpp$.
    Firstly, let $\sheaff$ be a sheaf on $\curve$ supported set theoretically on $\fppred$.
    Then we can find a short exact sequence:
    \[
        0 \to \sheaff_1 \to \sheaff \to \sheaff_2 \to 0
    \]
    where $\operatorname{supp} \sheaff_1$ is contained in $\{u = s^r = 0\}$ for some $r \geq 1$ and $\operatorname{supp} \sheaff_2 \subset \fpp$.
    Now assume that $\operatorname{supp} \sheaff$ is a thickening of $\fpp$, so that $\sheaff \not\cong \sheaff_2$ and hence $\sheaff_1 \neq 0$.
    Since $\operatorname{supp} \sheaff_1 \subset \curvebase$ we find that $\sigma_{*}\sheaff_1 \neq 0$.
    Now, left exactness of $\sigma_{*}$ implies that $\sigma_{*}\sheaff \neq 0$.
    Hence, any sheaf $\sheaff \in \ker \sigma_{*}$ has scheme-theoretic support contained in $\fpp$.

    Finally, if $\sheaff^{\bullet} \in \db{\curve}$ is a complex, we use a standard spectral sequence argument (see, e.g. \citep[Lemma 3.1]{bridgeland02}):
    \[E_2^{p, q} = R^{p}\sigma_{*}(H^{q}(\sheaff^{\bullet})) \implies R^{p + q}\sigma_{*}\sheaff^{\bullet}.\]
    This degenerates on the second page, so that if $\operatorname{supp} \sheaff^{\bullet}$
    is a thickening of $\fpp$, then the non-vanishing of $\sigma_{*}H^{q}(\sheaff^{\bullet})$ for some $q$ implies the non-vanishing of $\sigma_{*}\sheaff^{\bullet}$.
\end{proof}

Now we let $\sheaf(a, b)$ denote the restriction of the line bundle on $\pp \times \pp$ with degree $a$ along $\pp_{s:t}$ and degree $b$ along $\pp_{u:v}$.
Then, the object $\sh{P}(-1)$ has an infinite locally free resolution
\[
0 \leftarrow \sh{P}(-1) \leftarrow \sheaf(0, -1) \xleftarrow{\cdot s} \sheaf(-1, -1) \xleftarrow{\cdot su} \sheaf(-2, -2) \xleftarrow{\cdot s} \sheaf(-3, -2) \leftarrow \cdots
\]
so that $\Ext^{i}(\sh{P}(-1), \sh{P}(-1))$ can be computed using the chain complex:
\[
    0 \to \sheaf \xrightarrow{0} \sheaf \xrightarrow{0} \sheaf(1) \xrightarrow{0} \sheaf(1) \xrightarrow{0} \sheaf(2) \xrightarrow{0} \sheaf(2) \xrightarrow{0} \cdots
\]
where all maps are zero, since the coordinate $s$ becomes $0$ when restricted to $\fppred$.
Hence, $\mathcal{E}xt^{i}(\sh{P}(-1), \sh{P}(-1))$ is the pushforward of a line bundle on $\fppred$, and to compute global Ext we can simply take global sections as the local-to-global spectral sequence degenerates on the second page.
If $n$ is even, then this shows that $\Ext^{n}(\sh{P}(-1), \sh{P}(-1))$ has dimension $n/2 + 1$, with basis given by $\{ u^{i}v^{j} \;\vert\; i + j = n / 2\}$.
If $n$ is odd, then $\Ext^{n}(\sh{P}(-1), \sh{P}(-1))$ has dimension $(n + 1)/2$ with basis given by $\{ u^iv^j \;\vert\; i + j = (n - 1)/2\}$.
Concretely, the dimensions of $\Ext^{n}(\sh{P}(-1), \sh{P}(-1))$ form a sequence: $(1, 1, 2, 2, 3, 3, \dots)$.
We guess from this observation that the algebra structure is that of a polynomial ring with a generator of degree $1$ and a generator of degree $2$; this guess is verified in the next theorem.

\begin{defn}
    An object $P \in \db{X}$ for a variety $X$ is called a \term{compound} $\mathbb{P}^{\infty}$-object, or $\cpinf$, if there is an isomorphism of graded algebras:
    \[
    \Extcomplex(P, P) \cong \mathbb{C}[\epsilon, \theta]
    \]
    where $\deg \epsilon = 1$ and $\deg \theta = 2$.
\end{defn}

\begin{thm}\label{thm:curveabsorption}
    The object $P = \sh{P}(-1)$ is a $\cpinf$ object on $\curve$.
\end{thm}
\begin{proof}
    Consider the locally free resolution for $P$. The element of $\Ext^1(P, P) \cong \mathbb{C}$ arises from the mapping \[\mathbbm{1} : \sheaf(-1, -1) \to P\] and we lift this to a morphism $\epsilon : P \to P[1]$ by defining a chain map on the locally free resolution, so that all triangles and quadrilaterals in the diagram below commute:
    \[\begin{tikzcd}
	P && {\mathcal{O}(-1, -1)} & {\mathcal{O}(-2, -2)} & \cdots \\
	& {\mathcal{O}(0, -1)} & {\mathcal{O}(-1, -1)} & {\mathcal{O}(-2, -2)} & \cdots
	\arrow["{\cdot su}"', from=1-4, to=1-3]
	\arrow["{\cdot su}", from=2-4, to=2-3]
	\arrow["{\mathbbm{1}}"', from=1-3, to=1-1]
	\arrow["{\operatorname{res}}"', from=2-2, to=1-1]
	\arrow["{\cdot t}"{description}, from=1-3, to=2-2]
	\arrow["{\cdot ut}"{description}, from=1-4, to=2-3]
	\arrow["{\cdot s}", from=2-3, to=2-2]
	\arrow["{\cdot s}"', from=1-5, to=1-4]
	\arrow["{\cdot s}", from=2-5, to=2-4]
	\arrow["{\cdot t}"{description}, from=1-5, to=2-4]
\end{tikzcd}\]
Next, $\Ext^2(P, P) \cong \mathbb{C}^{\oplus 2}$ has basis given by the elements $\{u, v\}$.
We lift the element given by $v$ to a chain map $\theta : P \to P[2]$ similarly:
\[\begin{tikzcd}
	P && {\mathcal{O}(-2, -2)} & {\mathcal{O}(-3, -2)} & \cdots \\
	& {\mathcal{O}(0, -1)} & {\mathcal{O}(-1, -1)} & {\mathcal{O}(-2, -2)} & \cdots
	\arrow["{\cdot s}"', from=1-4, to=1-3]
	\arrow["{\cdot su}", from=2-4, to=2-3]
	\arrow["{\cdot v}"', from=1-3, to=1-1]
	\arrow["{\operatorname{res}}"', from=2-2, to=1-1]
	\arrow["{\cdot vt^2}"{description}, from=1-3, to=2-2]
	\arrow["{\cdot vt^2}"{description}, from=1-4, to=2-3]
	\arrow["{\cdot s}", from=2-3, to=2-2]
	\arrow["{\cdot su}"', from=1-5, to=1-4]
	\arrow["{\cdot s}", from=2-5, to=2-4]
	\arrow["{\cdot vt^2}"{description}, from=1-5, to=2-4]
\end{tikzcd}\]
One now sees that $\epsilon \circ \theta = \theta \circ \epsilon$ since both compositions are given by the chain map lifting the element $v \in \Ext^3(P, P)$:
\[\begin{tikzcd}
	P && {\mathcal{O}(-3, -2)} & {\mathcal{O}(-4, -3)} & \cdots \\
	& {\mathcal{O}(0, -1)} & {\mathcal{O}(-1, -1)} & {\mathcal{O}(-2, -2)} & \cdots
	\arrow["{\cdot su}"', from=1-4, to=1-3]
	\arrow["{\cdot su}", from=2-4, to=2-3]
	\arrow["{\cdot v}"', from=1-3, to=1-1]
	\arrow["{\operatorname{res}}"', from=2-2, to=1-1]
	\arrow["{\cdot vt^3}"{description}, from=1-3, to=2-2]
	\arrow["{\cdot uvt^3}"{description}, from=1-4, to=2-3]
	\arrow["{\cdot s}", from=2-3, to=2-2]
	\arrow["{\cdot s}"', from=1-5, to=1-4]
	\arrow["{\cdot s}", from=2-5, to=2-4]
	\arrow["{\cdot vt^3}"{description}, from=1-5, to=2-4]
\end{tikzcd}\]
In general, we claim that the lift of the element $u^{i}v^{j} \in \Ext^{n}(P, P)$ is given by the chain map $\epsilon^{2i}\theta^{j}$ if $n$ is even, and $\epsilon^{2i + 1}\theta^{j}$ if $n$ is odd.
This is verified by observing that in the $n$ even case, the chain map lifting $\epsilon^{2i}\theta^j$ is given by:
\[\begin{tikzcd}
	P & {\mathcal{O}(-n, -(n/2)-1)} & {\mathcal{O}(-n-1, -(n/2)-1)} & \cdots \\
	& {\mathcal{O}(2j - n, j - (n/2) - 1)} & {\mathcal{O}(2j - n - 1, j - (n/2) -1)} & \cdots \\
	P & {\mathcal{O}(0, -1)} & {\mathcal{O}(-1, -1)} & \cdots
	\arrow["{{\cdot s}}"', from=2-3, to=2-2]
	\arrow["{{\cdot s}}"', from=3-3, to=3-2]
	\arrow["{{\cdot su}}"', from=2-4, to=2-3]
	\arrow["{{{\cdot v^jt^{2j}}}}", from=1-2, to=2-2]
	\arrow["{\cdot s}"', from=1-3, to=1-2]
	\arrow["{{\cdot su}}"', from=1-4, to=1-3]
	\arrow["{{{\cdot v^jt^{2j}}}}", from=1-3, to=2-3]
	\arrow["{\cdot su}"', from=3-4, to=3-3]
	\arrow["{{\cdot u^it^{2i}}}", from=2-2, to=3-2]
	\arrow["{{\cdot u^iv^j}}"', from=1-2, to=1-1]
	\arrow["{{\operatorname{res}}}"{description}, from=3-2, to=3-1]
	\arrow["{{\cdot u^it^{2i}}}", from=2-3, to=3-3]
	\arrow[Rightarrow, no head, from=1-1, to=3-1]
\end{tikzcd}\]
while in the odd case, the chain map lifting $\epsilon^{2i + 1}\theta^{j}$ is given by the following.
\[\begin{tikzcd}
	P & {\mathcal{O}(-n - 1, -(n + 1)/2)} & {\mathcal{O}(-n-2, -(n+ 3)/2)} & \cdots \\
	& {\mathcal{O}(2j - n - 1, j - (n + 1)/2)} & {\mathcal{O}(2j - n - 2, j - (n + 3)/2)} & \cdots \\
	P & {\mathcal{O}(0, -1)} & {\mathcal{O}(-1, -1)} & \cdots
	\arrow["{{\cdot s}}"', from=3-3, to=3-2]
	\arrow["{{\cdot su}}"', from=3-4, to=3-3]
	\arrow["{{\operatorname{res}}}"{description}, from=3-2, to=3-1]
	\arrow["{{\cdot u^iv^j}}"', from=1-2, to=1-1]
	\arrow["{{{\cdot v^jt^{2j}}}}", from=1-2, to=2-2]
	\arrow["{{\cdot u^it^{2i + 1}}}", from=2-2, to=3-2]
	\arrow[Rightarrow, no head, from=1-1, to=3-1]
	\arrow["{\cdot su}"', from=1-3, to=1-2]
	\arrow["{{{\cdot v^jt^{2j}}}}", from=1-3, to=2-3]
	\arrow["{{\cdot u^{i + 1}t^{2i + 1}}}", from=2-3, to=3-3]
	\arrow["{{\cdot s}}"', from=1-4, to=1-3]
	\arrow["{{\cdot su}}"', from=2-3, to=2-2]
	\arrow["{{\cdot s}}"', from=2-4, to=2-3]
\end{tikzcd}\]

\end{proof}

\begin{remark}\label{remark:formal}
    In fact, the above proof shows that the algebra $\operatorname{RHom}(P, P)$ is formal, i.e. quasi-isomorphic to $\Extcomplex(P, P)$; the chain maps constructed above give us the quasi-isomorphism.
\end{remark}

\begin{remark}
    There is in fact a general class of examples where one can perform the same trick, given by considering the trivially thickened curve $\pp \times \mathbb{C}[s]/s^n$ intersecting a smooth curve transversely. 
    For $n > 2$, with $P$ defined analogously, the algebra $\operatorname{RHom}(P, P)$ will not be formal.
\end{remark}

Note that $\omega_{\curve} \cong \sheaf(0, -1)$.
Hence, by the Serre duality statement of \Cref{serreduality}, we also have the following semiorthogonal decomposition for $\curve$, which we will use in the sequel.
\begin{equation} \label{eq:sodforC}
    \db{\curve} = \sod{\db{\curvebase}, \sh{P}(-1) \otimes \omega_\curve^{-1}} =
\sod{\db{\curvebase}, \sh{P}}.
\end{equation}
Note here that $\sh{P}$ is the skyscraper sheaf along $\fppred$, and we will define
\begin{equation}\label{eq:absorbingcat}
    \sh{D} := \langle \sh{P} \rangle \subset \db{X}
\end{equation}
to be the absorbing subcategory of $\db{X}$.

\begin{remark}
Given a $\cpinf$ object $P$, we expect that there is an equivalence given by Koszul duality
\begin{align*}
    & \gens{P} \simeq  \mathbf{D}^{\operatorname{perf}}(\mathbb{C}[\epsilon, \theta]) \simeq \db{\mathbb{C}[w, r]/(w^2, r^2)} && \deg \epsilon = 1, \deg \theta = 2, \\ &&& \deg w = 0, \deg r = -1.
\end{align*}
This would follow by computing $\operatorname{Ext}^{\bullet}(\mathbb{C}, \mathbb{C})$ over the algebra $\mathbb{C}[w, r]/(w^2, r^2)$.
After verifying that the dg-algebras involved are formal, for example using a similar argument to the one given in \Cref{remark:formal}, the desired Koszul duality statement would then follow from \cite[Theorem 8.5(b)]{keller}.
\end{remark}

\section{Semiorthogonal decomposition for a resolution of singularities} \label{section:catres}

Let $B$ be a smooth threefold which contains the curve $\curve$ considered in \Cref{section:cpinfobjects} as a locally complete intersection.
In this section, to allow us to perform explicit calculations, we will also assume that $B$ contains an open subset 
\begin{equation}\label{eq:defub}
    U_B \cong \mathbb{A}^2_{x, s} \times \pp_{u:v}
\end{equation}
such that $\curve \cap U_B$ is just the curve $\curve$ missing a point at infinity on the reduced branch.
For example, we may take \[B = \pp_{x:y} \times \pp_{s:t} \times \pp_{u:v}\] which contains $\curve$ as the vanishing locus $V(x, s^2u)$ and $U_B$ as the locus $\{y \neq 0, t \neq 0\}$.
Let $X$ be the threefold constructed through the blow-up:
\[
X := \bll{\curve}B.
\]
Then, $X$ is singular with a line of surface nodes compounded with a threefold nodal singularity at one point. 
In particular, this singularity is non-isolated.
By Orlov's blow-up formula given in \cref{eq:sodblowup} combined with the absorption for $\curve$ given by \cref{eq:sodforC}, we immediately obtain a semiorthogonal decomposition:
\[
\db{X} = \sod{\db{\curve}, \db{B}} = \sod{\db{\pp_{s:t}}, \sh{D}, \db{B}}
\]
where $\sh{D}$ is the absorbing category as in \eqref{eq:absorbingcat}, generated by the $\cpinf$ object $\sh{P}$.

It is worth considering how each component of the decomposition is embedded, using \eqref{eq:blowupembedding}.
Firstly, we define a line bundle on $\curve$
\[
\sh{L}(a) := \sheaf(a, 0)
\]
which is the pullback of a line bundle $\sheaf(a)$ via the contraction map $\curve \to \curvebase \cong \pp_{s:t}$.
Note $\sh{L}(a)$ has degree $a$ along the $\pp_{s:t}$ branch and is trivial along the thickened branch.

Next, $\db{\curve}$ is embedded into $X$ by the functor $\Phi: \db{\curve} \to \db{X}$ which pulls back to the exceptional locus, tensors with the $\sheaf(-1)$ bundle and pushes forward into $X$.
Since $X$ is the blow-up of $B$ along $\curve$, we observe that the exceptional locus has two components, both of which are $\pp$-bundles.
\begin{align}
    & F \to \pp_{s:t} \label{def:F}\\
    & G \to \pp_{u:v} \times \fatpoint\label{def:G}
\end{align}
We can also pullback along $\pp \xhookrightarrow{} \pp \times \fatpoint$ to obtain 
\begin{equation}
    G_{\red} \to \pp_{u:v}.\label{def:Gred}
\end{equation}
Then, we can write our semiorthogonal decomposition as
\begin{equation}\label{eq:sodx}
    \db{X} = \sod{\Phi(\sh{L}(1)), \Phi(\sh{L}(2)), \Phi(\sh{P}), \db{B}}
\end{equation}
where the sheaves $\Phi(\sh{L}(a))$ are supported on $F \cup G$ and $\Phi(\sh{P})$ is supported on $G_{\red}$.
Note that we have chosen $\db{\pp_
{s:t}} = \sod{\sheaf(1), \sheaf(2)}$.

Locally over $U_B$, $X$ takes the form
\[U_X := \{xb = s^2ua\} \subset \mathbb{P}_{U_b}(\sheaf(1)_a \oplus \sheaf_b)\]
and hence has a compound $A_2$ singularity of the form described previously.
We find that our semiorthogonal decomposition provides a categorical absorption for this singularity.
The goal of this section is to show that the semiorthogonal decomposition descends from one for a resolution for $X$.

\begin{remark}
    If we instead start with $\pp \times \mathbb{C}[s]/s^n$ intersecting a $\pp$ transversely, then we obtain a threefold with a non-isolated compound $A_n$ singularity locally given by the equation
    \[
    \{xb = s^nu\} \subset \mathbb{A}^4
    \]
    so that we have a line of surface $A_{n - 1}$ singularities, compounded with a threefold node at the origin. 
\end{remark}

\subsection{Constructing a geometric resolution}

The process of constructing the resolution is as follows.
To resolve the singularities, we first blow-up along the exceptional component $F$ defined in \eqref{def:F}:
\[
\widehat{f}_X : \widehat{X} \to X
\]
to obtain a singular threefold $\widehat{X}$.
Since $F$ intersects the singular locus only at the threefold node, we find that $\widehat{X}$ has a $\pp$-family of surface node singularities.
To resolve these, we blow-up along the strict transform of $G_{\red}$ defined in \eqref{def:Gred}:
\[
\varphi_{\widehat{X}} : \widetilde{X} \to \widehat{X}
\]
to obtain a smooth threefold $\widetilde{X}$.
Define
\begin{equation}\label{def:sigma}
    \sigma : \widetilde{X} \to X
\end{equation}
to be the composition $\widetilde{X} \to \widehat{X} \to X$.

It may be helpful to refer to the surface analogue: we may consider the blow-up
\[
\bl{x^2, y} \mathbb{A}^2_{x, y} \cong \{ x^2a = yb \} \subset \mathbb{A}^2 \times \pp_{a:b}
\]
so that the exceptional divisor $V(x^2, y)$ is a thickened $\pp$, and the space has a nodal singularity at the point $(0, 0, [0:1])$.
Then blowing up the reduced exceptional divisor $V(x, y)$ gives a minimal resolution of the singularity; our second blow-up $\widetilde{X} \to \widehat{X}$ is a family version of this process.

\begin{figure}[!ht]
    \centering
    \includegraphics[page=1,width=0.8\textwidth]{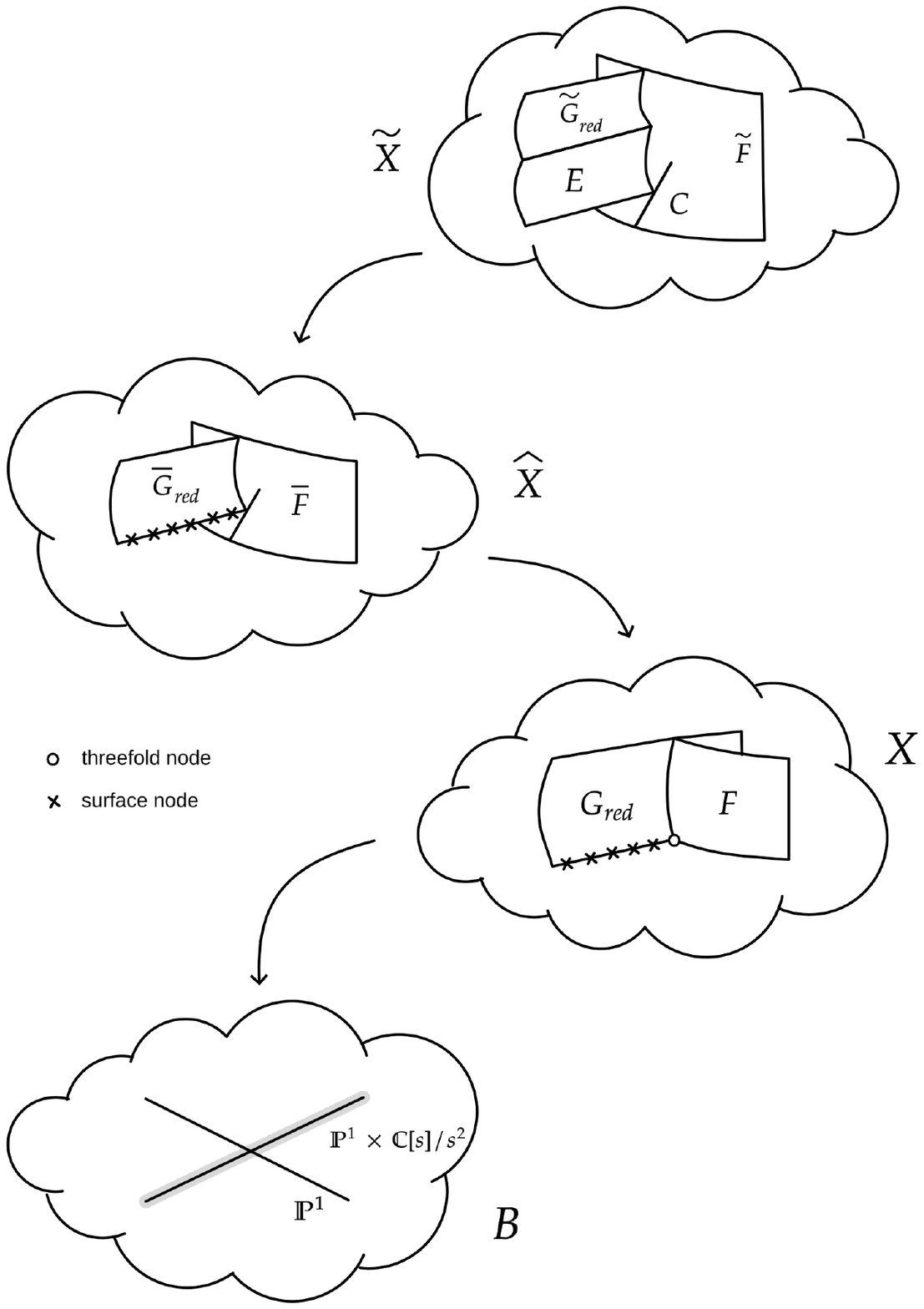}
    \caption{Constructing the resolution. Here, $\overline{F}, \overline{G}_{\red}$ denote the strict transforms in $\widehat{X}$ and $\widetilde{F}, \widetilde{G}_{\red}$ denote the strict transforms in $\widetilde{X}$.}
    \label{fig:resolution-construction}
\end{figure}

\subsubsection{Local analysis}
We illustrate this construction locally by considering the open subset $U_B$ defined in \eqref{eq:defub}.
Recall that we obtain an open subset of $X$:
\[U_X := \{xb = s^2ua\} \subset \mathbb{P}_{U_b}(\sheaf(1)_a \oplus \sheaf_b) =: P_X.\]
Locally, $F \cap U_X \cong V(x, u)$ and $G \cap U_X \cong V(s^2, x)$.
To obtain an open in $\widehat{X}$, we blow-up at the ideal $\langle x, u \rangle$:
\[
U_{\widehat{X}} := \left\{\begin{array}{c}bd = s^2ac\\xc = ud\end{array}\right\} \subset 
    \pbundle_{P_{X}}( \sheaf_c \oplus \base{1}_d ) =: P_{\widehat{X}}.
\]
where we denote $\pi_B : P_X \to U_B$.
We can now resolve the singularities by blowing up at $\langle d, s \rangle$.
\[
U_{\widetilde{X}} := \left\{\begin{array}{c}
    bf = sace\\
    xc = ud\\
    de = sf
    \end{array}\right\} \subset
        \pbundle_{P_{\widehat{X}}}(\sheaf(1)_e \oplus \base{1}_f) =: P_{\widetilde{X}}
\]
where this time we are denoting $\pi_B : P_{\widehat{X}} \to U_B$ and $\sheaf(1)$ is the tautological line bundle on $P_{\widehat{X}}$.
As a GIT quotient, $P_{\widetilde{X}}$ is constructed as the quotient $\mathbb{C}^{10}/(\mathbb{C}^{*})^4$ with weights shown in \Cref{tab:gitquotient}.

\begin{table}
    \centering
    \begin{tabular}{cccccccccc}
         x&  s&  u&  v&  a&  b&  c&  d&  e& f\\
         \hline
         0&  0&  1&  1&  -1&  0&  0&  -1&  0& -1\\
         &  &  &  &  1&  1&  0&  0&  0& 0\\
         &  &  &  &  &  &  1&  1&  -1& 0\\
         &  &  &  &  &  &  &  &  1& 1\\
    \end{tabular}
    \caption{Matrix of weights for $P_{\widetilde{X}}$, with omitted values equal to $0$}
    \label{tab:gitquotient}
\end{table}

Note that the exceptional locus of $U_{\widehat{X}} \to U_X$ is a curve $V(x, u, s, b) \cong \pp_{c:d}$.
We denote by 
\begin{equation}
    C = V(x, u ,s, b, e) \cong \pp_{c:d} \subset U_{\widetilde{X}}
\end{equation}
the strict transform of this curve in $\widetilde{X}$.
The strict transform $\overline{G}_{\red}$ of $G_{\red}$ in $\widehat{X}$ intersects the singular locus along the section where the fibre coordinate $b$ vanishes.
The effect of blowing up this locus is to introduce an exceptional surface $E \subset U_{\widetilde{X}}$.
The exceptional locus of the composition $\sigma : \widetilde{X} \to X$ is the union $E \cup C$, and it is precisely the preimage of the reduced singular locus $V(x, b, s) \subset X$.
Note that the two subvarieties $E$ and $C$ intersect transversely at a point.

\begin{lemma}\label{lemma:exceptionalhirzebruch}
    Let $\widetilde{G}_{\red}$ be the strict transform of $G_{\red}$ in $\widetilde{X}$.
    Then we have
    \begin{align*}
        \widetilde{G}_{\red} & \cong \pbundle_{\pp_{u:v}}(\sheaf_a \oplus \sheaf(1)_b) \\
        E & \cong \pbundle_{\pp_{u:v}}(\sheaf_e \oplus \sheaf(1)_f).
    \end{align*}
\end{lemma}
\begin{proof}
    In $\widehat{X}$, the locus $\overline{G}_{\red}$ is given by $\overline{G}_{\red} = V(x, d, s) \subset U_{\widehat{X}}$, and so its strict transform is given by $\widetilde{G}_{\red} = V(x, d, s, f) \subset U_{\widetilde{X}}$. 
    The second statement follows from similarly noting that $E = V(x, d, s, b) \subset U_{\widetilde{X}}$.
\end{proof}

We also see that the strict transform $\overline{F} \subset \widehat{X}$ of $F$ intersected with $U_{\widehat{X}}$ is a nodal surface:
\[
\{ bd = s^2ac \} \subset \mathbb{A}^1_s \times \pp_{a:b} \times \pp_{c:d}
\]
which is isomorphic to the blow-up of $\mathbb{A}^1_s \times \pp_{a:b}$ at a fattened point cut out by the ideal $\langle s^2a, b \rangle$.
Denoting by $\widetilde{F}$ the strict transform of $F$ in $\widetilde{X}$, we find that $\widetilde{F}$ is a ruled surface blown-up at a point, and then blown-up again at a point on the exceptional curve.

\subsubsection{Alternative construction of the resolution}
A key computational tool will be the following alternative construction of $\widetilde{X}$, where we only blow-up smooth varieties along a $\pp$.
This allows us to easily obtain a concrete semiorthogonal decomposition for $\db{\widetilde{X}}$.

\begin{figure}[!ht]
    \centering
    \includegraphics[page=2,width=0.8\textwidth]{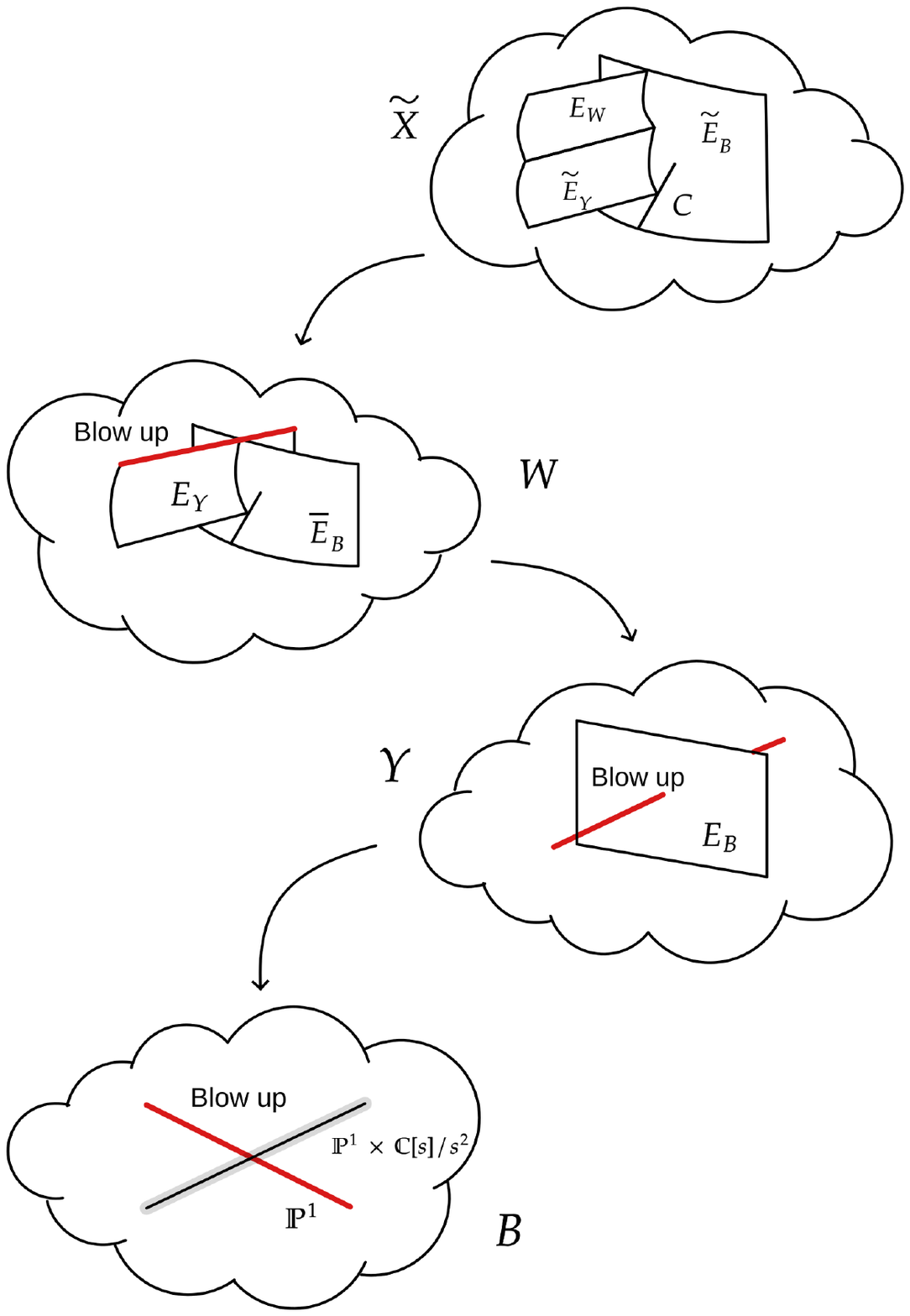}
    \caption{Alternative construction of the resolution, with blow-up centers indicated in red. Here, $\overline{E}_B$ denotes the strict transform of $E_B$ to $W$, and $\widetilde{E}_B, \widetilde{E}_Y$ denote the strict transforms of $E_B, E_Y$ respectively to $\widetilde{X}$.}
    \label{fig:alt-construct}
\end{figure}

\begin{prop}\label{prop:altres}
    There are varieties $W, Y$ and a sequence of contractions $\varphi_W, \psi, h_Y$ where each map is a blow-up of a smooth rational curve, such that the following diagram commutes.
    \begin{equation}\label{eq:resolutiondiagram}
    \begin{tikzcd}
	{\widetilde{X}} && W \\
	\\
	{\widehat{X}} && Y \\
	\\
	X && B
	\arrow["{\varphi_{\widehat{X}}}"', from=1-1, to=3-1]
	\arrow["{\widehat{f}_X}"', from=3-1, to=5-1]
	\arrow["{h_Y}", from=3-3, to=5-3]
	\arrow["{h_X}"', from=5-1, to=5-3]
	\arrow["{\widehat{f}_Y}", from=3-1, to=3-3]
	\arrow["\psi", from=1-3, to=3-3]
	\arrow["{\varphi_W}", from=1-1, to=1-3]
	\arrow["{\varphi_Y}", from=1-1, to=3-3]
    \end{tikzcd}
    \end{equation}
    Furthermore, letting $E_B, E_Y, E_W$ be the exceptional loci of $h_Y, \psi, \varphi_W$ respectively and $\widetilde{E}_B, \widetilde{E}_Y$ the strict transforms of $E_B, E_Y$ respectively to $\widetilde{X}$, we have that:
    \begin{enumerate}
        \item $\widetilde{E}_Y$ coincides with $E$;
        \item $E_W$ coincides with $\widetilde{G}_{\red}$;
        \item $\widetilde{E}_B$ coincides with $\widetilde{F}$.
    \end{enumerate}
\end{prop}
\begin{proof}
    The map $h_Y : Y \to B$ is obtained by blowing-up $B$ along the \textit{reduced} branch \[\curvebase \cong \pp_{s:t} \subset \curve \subset B\] to obtain a variety $h_Y : Y \to B$.
    Now, $\widehat{X}$ can be obtained by blowing up the strict transform $\overline{\fpp}$ of the non-reduced branch $\fpp \subset \curve \subset B$.
Hence we have a map $\varphi_Y : \widetilde{X} \to Y$, and we see that the total transform of the ideal cutting out the reduced curve $\overline{\fppred}$ is locally principal in $\widetilde{X}$.
Letting $\psi : W \to Y$ be the blow-up in $\overline{\fppred}$, we see that $\widetilde{X} \to Y$ factors as $\widetilde{X} \to W \to Y$ by the universal property of blow-ups.

Note that outside of the open subset $U_B$ both constructions of $\widetilde{X}$ only blow-up $B$ at the single point at infinity on $\curvebase \subset \curve$, and otherwise all blow-ups have centers contained in some open over $U_B$.
Define the open in $Y$ over $U_B$:
\[
U_Y := \{ xc = ud \} \subset \pbundle_{U_B}(\sheaf_c \oplus \sheaf(1)_d) =: P_{Y}.
\]
Similarly, there is an open subset of $W$ found by blowing up $P_Y$ at $\gens{d, s}$ and pulling back:
\[
U_W := \left\{\begin{array}{c}xc = ud\\de = sf\end{array}\right\} \subset 
    \pbundle_{P_Y}(\sheaf(1)_e \oplus \base{1}_f) =: P_W.
\]
where we denote $\pi_B : P_Y \to U_B$.
Finally, we blow up $P_W$ along the ideal $\gens{f, sce}$ and pull back the equations for $U_W$:
\[
U_{\widetilde{X}} = \left\{\begin{array}{c}xc = ud\\de = sf\\bf = (sce)a\end{array}\right\} \subset \pbundle_{P_W}(\base{1}_a \oplus \sheaf_b) =: P_{\widetilde{X}}^{\prime}
\]
where we denote $\pi_B : P_W \to U_B$.
It is clear that $P^{\prime}_{\widetilde{X}} \cong P_{\widetilde{X}}$ (e.g. by writing out the matrix of weights for $P^{\prime}_{\widetilde{X}}$ and comparing with \Cref{tab:gitquotient}), and so we recover $U_{\widetilde{X}}$.

The center of $\psi : W \to Y$ is given by $Z_W = V(x, d, s) \subset U_Y$, and the exceptional locus of $\psi : W \to Y$ is given by
\[
E_Y = V(x, d, s) \cong \pbundle_{\pp_{u:v}}(\sheaf_e \oplus \sheaf(1)_f) \subset U_W.
\]
The center of the blow up $\widetilde{X} \to W$ is given by $Z_{\widetilde{X}} = V(f, d, x, s) \subset U_W$ which is a copy of $\pp_{u:v}$ arising as a section of $E_Y \to \pp_{u:v}$, and we see that 
\[
E_W = V(f, d, s, x) \subset U_{\widetilde{X}}
\]
is the exceptional divisor of the blow-up $\widetilde{X} \to W$.
Now, it is clear by considering the equations cutting out each subvariety that $E_W$ coincides with $\widetilde{G}_{\red}$ and $\widetilde{E}_Y$ coincides with $E$ in $\widetilde{X}$.
It is also clear that $\widetilde{E}_B$ coincides with $\widetilde{F}$ in $\widetilde{X}$ since the both give the preimage of $\curvebase$ along the map $\widetilde{X} \to B$.
\end{proof}
The alternative construction is visualised in \Cref{fig:alt-construct}.

\subsubsection{Exceptional loci} \label{subsec:exceptionalloci}
Next, we consider the irreducible component $E$ of the exceptional locus of $\sigma$, recalling the notation for line bundles on a ruled surface introduced in \cref{eq:ruledpicard}.
\begin{lemma}\label{restrictcanonical}
    For the divisor $E \subset \widetilde{X}$, we have that:
    \begin{align*}
        & \sh{N}_{E / \widetilde{X}} \cong \sheaf(-2, 1) \nonumber \\
        & \omega_{\widetilde{X}}\restrict{E} \cong \sheaf(0, -2)
    \end{align*}
\end{lemma}
\begin{proof}
    Recall that $E \cong \mathbb{P}_{\mathbb{P}^{1}_{u:v}}(\sheaf_{E} \oplus \sheaf(1)_f)$ by \Cref{lemma:exceptionalhirzebruch}.
    The bundle $\sh{N}_{E/\widetilde{X}}$ has degree $-2$ along the fibres of the bundle $E \to \pp_{u:v}$, which we can see by considering that each fibre is the exceptional curve in the resolution of a surface node.
    Also note that since $E \subset \widetilde{X}$ is the strict transform of $E_Y \subset W$, $E$ is the blow up of $E_Y$ along the section given by $V(f) \subset E_Y$. 
    This implies that the restriction of $\sh{N}_{E/\widetilde{X}}$ to the section $V(e)$ is the same as the restriction of $\sh{N}_{E_Y/W} \cong \sheaf(-1)$ to $V(e)$.
    Hence, we find that $\sh{N}_{E/\widetilde{X}} \cong \sheaf(-2, 1)$, and using the formula for $\omega_{E}$ given in \eqref{eq:ruledcanonical} we compute that:
\begin{align*}
    \omega_{\widetilde{X}}\restrict{E} 
    & \cong \omega_{E} \otimes \sh{N}^{\vee}_{E/\widetilde{X}} \\ 
    & \cong \sheaf(-2, -1) \otimes \sheaf(2, -1) \cong \sheaf(0, -2).
\end{align*}
\end{proof}

Finally, we consider the curve $C$.
Recall the definition of a spherical object from \Cref{sec:spherical}.
\begin{lemma}\label{lemma:threespherical}
    The object $\sheaf_{C}(-1)$ is a $3$-spherical object in $\db{\widetilde{X}}$.
\end{lemma}
\begin{proof}
    This follows if we can contract $C$ to an ordinary double point.
    Starting from $B$, let $\curve_{\operatorname{node}}$ be the union of the reduction of the branches, and consider the blow-up:
    \[
        Y' = \bll{\curve_{\operatorname{node}}}B
    \]
    so that $Y'$ is a nodal variety.
    Now note that the total transform of the locus $\curve_{\operatorname{node}}$ in $Y$ is the union of $E_B$, which is a Cartier divisor, and the center $Z_W$ of the blow-up $W \to Y$.
    It follows from the universal property of the blow-up that $W$ can be written as a blow-up of $Y'$ at the strict transform of $\curvebase$, which is a ruled surface in $Y'$ intersecting the node.
    In particular, there is an exceptional curve $C_1$ which contracts to give the node.
    The blow-up $\varphi_W : \widetilde{X} \to W$ has center disjoint from $C_1$, and we can compute by working locally over $U_B$ that $C_1$ is identified with the curve $C$ under the blow-up. 
    Hence, $C$ can also be contracted to a threefold node.
\end{proof}

\subsection{A full exceptional sequence on the resolution}

With the geometric setup established, we may now analyse $\db{\widetilde{X}}$ using Orlov's blow-up formula once again.
We recall the exceptional loci in the alternative construction from \Cref{subsec:exceptionalloci}, as well as the maps described in diagram \eqref{eq:resolutiondiagram}.
The center of each blow-up in the alternative construction is a locally complete intersection. 
This gives us the following semiorthogonal decompositions, where for a $\pp$-bundle $\pi : \Sigma \to \pp$ embedded into a variety $T$ as $j : \Sigma \to T$ we denote by $\sheaf_{\Sigma}(a, b)$ the pushforward $j_{*}(\sheaf(a) \otimes \pi^{*}\sheaf(b))$.
\begin{align*}
    \db{Y} & = \sod{ \db{\pp_{s:t}}, \db{B} } \\
        & = \sod{ \sheaf_{E_B}(-1, -1), \sheaf_{E_B}(-1, 0), h_Y^{*}\db{B} }
\end{align*}
\begin{align*}
    \db{W} & = \sod{ \db{\pp_{u:v}}, \db{Y} } \\
    & = \sod{ \sheaf_{E_Y}(-1, -1), \sheaf_{E_Y}(-1, 0), \psi^*\db{Y} }
\end{align*}
\begin{align*}
    \db{\widetilde{X}} & = \sod{ \db{\pp_{u:v}}, \db{W} } \\
    & = \sod{ \sheaf_{E_W}(-1, -1), \sheaf_{E_W}(-1, 0), \varphi_{W}^*\db{W} }
\end{align*}

The component in the semiorthogonal decomposition equivalent to $\db{B}$ is immediately seen to be orthogonal to $\ker \sigma_{*}$.
We therefore consider the following objects
\begin{align}\label{def:exceptionals}
    & \exobj_{0} := \sheaf_{E_W}(-1, -1) && \exobj_{1} := \sheaf_{E_W}(-1, 0) \nonumber \\
    & \sheaff_{0} := \varphi_W^{*}\sheaf_{E_Y}(-1, -1) && \sheaff_{1} := \varphi_W^{*}\sheaf_{E_Y}(-1, 0) \\
    & \sh{A}_0 := \varphi_Y^{*}\sheaf_{E_B}(-1, -1) && \sh{A}_1 := \varphi_Y^{*}\sheaf_{E_B}(-1, 0) \nonumber
\end{align}
so that $\ker \sigma_{*}$ lies in the exceptional sequence
\[
\ker \sigma_{*} \subset \sod{\exobj_0, \exobj_1, \sheaff_0, \sheaff_1, \sh{A}_0, \sh{A}_1}.
\]

It is useful to obtain a geometric description for the objects $\sh{F}_i$ and $\sh{A}_i$.
The following lemmas are certainly specialisations of more general results on the properties of pulling back skyscraper sheaves along blow-ups, but we give the concrete calculations here for completeness, using the properties of ruled surfaces given in \Cref{sec:ruledsurfaces}.

\begin{lemma}
    For $i \in \{0, 1\}$, $\sh{F}_i$ is the pushforward of the line bundle on $E \cup E_W$ which is:
    \begin{enumerate}
        \item The line bundle $\sheaf(-1, i - 1)$ on $E$
        \item The line bundle $\sheaf(0, i - 1)$ on $E_W$.
    \end{enumerate}
\end{lemma}
\begin{proof}
    We show the statement for $\sh{F}_1$, using its description from \cref{def:exceptionals}, and the corresponding statement for $\sh{F}_0$ follows by twisting appropriately. 
    Firstly, taking the exact sequence on $W$
\[
0 \to \sheaf \to \sheaf(E_Y) \to \sheaf_{E_Y}(-1, 0) \to 0
\]
and applying the pullback $\varphi_{W}^{*}$ shows that we have an exact sequence
\[
0 \to \sheaf \to \sheaf(E + E_W) \to \varphi_{W}^{*}\sheaf_{E_Y}(-1, 0) \to 0
\]
on $\widetilde{X}$.
Note that this holds even though the functor $\varphi_W^{*}$ is derived, i.e. $\textbf{L}^{i}\varphi_W^{*}(\sheaf_{E_Y}(-1, 0))$ vanishes for $i \neq 0$.
We find that the objects $\sheaff_i$ are supported on $E \cup E_W$.
Note that $E \cap E_W = V(f) \subset E$, and $f$ is a section of $\sheaf(1, -1)$.
Hence 
\[\sheaf(E + E_W)\restrict{E} \cong \sheaf(-2, 1) \otimes \sheaf(1, -1) \cong \sheaf(-1, 0)\]
and since $E \cap E_W = V(b) \subset E_W$, where $b$ is a section of $\sheaf(1, 0)$, we have
\[
\sheaf(E + E_W)\restrict{E_W} \cong \sheaf(-1, 0) \otimes \sheaf(1, 0) \cong \sheaf.
\]
We hence find that $\sh{F}_1$ is a pushforward of the line bundle on $E \cup E_W$ which is given by $\sheaf(-1, 0)$ on the $E$ component and $\sheaf$ on the $E_W$ component.
Note that the line bundle $\sheaf(a, b)$ on $E$ restricts to $\sheaf(b)$ on $E \cap E_W$, so such a line bundle exists on $E \cup E_W$.
\end{proof}

Recall that by \Cref{prop:altres}, $\widetilde{F} \subset \widetilde{X}$ coincides with the strict transform of $E_B$ in $\widetilde{X}$.

\begin{lemma}
    For $i \in \{0, 1\}$, the objects $\sh{A}_i$ are given as:
    \begin{align*}
    & \sh{A}_i \cong \sheaf_{\widetilde{F}}(-1, i - 1)
\end{align*}
\end{lemma}
\begin{proof}

For the object $\sh{A}_1$, we have a short exact sequence:
\[
0 \to \sheaf \to \sheaf(\widetilde{F}) \to \varphi_{Y}^{*}\sh{A}_1 \to 0
\]
so that $\sh{A}_1$ is the pushforward of the normal bundle of $\widetilde{F}$.
Recall that $\widetilde{F}$ is a ruled surface blown up twice, and has a chain of exceptional curves $C_E + C_W$, where we denote $C_E := \widetilde{F} \cap E$ and $C_W = \widetilde{F} \cap E_W$.
We denote by $\sheaf(z, w) := p^{*}\sheaf(z, w)$, where $p : \widetilde{F} \to E_B$ is the blow-up map.
Then the normal bundle must be of the form
\[
\sh{N}_{\widetilde{F} / \widetilde{X}} \cong \sheaf(z, w) \otimes \sheaf(aC_E + b C_W).
\]
We know that $z = -1$ and $w = 0$ since $\sh{N}_{E_B / Y} \cong \sheaf(-1, 0) $ and the restriction to the fibres of $\widetilde{F} \to E_B \to \pp_{u:v}$ is generically unchanged.
Next, one computes via the inclusion $C_E \xhookrightarrow{} E \xhookrightarrow{} \widetilde{X}$ that $\sh{N}_{C_E/ \widetilde{X}} \cong \sheaf \oplus \sheaf(-2)$.
Since $\sh{N}_{C_E/\widetilde{F}} \cong \sheaf(-2)$, we find that 
\[\sheaf(aC_E + bC_W)\restrict{C_E} \cong \sheaf \implies -2a + b = 0.\]
Similarly, $\sh{N}_{C_W/\widetilde{F}} \cong \sheaf(-1)$ and $\sh{N}_{C_W/\widetilde{X}} \cong \sheaf \oplus \sheaf(-1)$, so we find that
\[
\sheaf(aC_E + bC_W)\restrict{C_W} \cong \sheaf \implies -a + b = 0.
\]
Hence $b = a = 0$, and we obtain the statement for $\sh{A}_1$, with the case of $\sh{A}_0$ being similar.
\end{proof}

Observe now that na\"{i}vely pushing down this semiorthogonal decomposition via $\sigma_{*}$ does not recover the semiorthogonal decomposition we have for $X$.
Recall that we expressed our semiorthogonal decomposition as
\begin{equation*}
    \db{X} = \sod{\Phi(\sh{L}(1)), \Phi(\sh{L}(2)), \Phi(\sh{P}), \db{B}} 
\end{equation*}
where $\Phi : \db{\curve} \to \db{X}$ is the embedding induced by Orlov's blow-up formula.
We recall that $\sh{L}(a)$ denotes the line bundle $\sheaf(a, 0)$ on $\curve$ which has degree $a$ over the reduced branch and is trivial over the non-reduced branch, and $\Phi(\sh{P}) =: \sh{D}$ is the absorbing category.
We can see that $\sigma(E_W) = \sigma(E + E_W) = G_{\red} \subset X$, while $\sigma(\widetilde{F}) = F \subset X$.
Hence we immediately see that $\sigma_{*}\sh{E}_1$ is our $\cpinf$ object $\Phi(\sh{P})$, but $\sigma_{*}\sh{E}_0$ is a different $\cpinf$ object, namely $\Phi(\sh{P}(-1))$.
Meanwhile, the objects $\sh{A}_i$ are supported on $\widetilde{F}$ and hence push down to sheaves supported only on $F$.

It turns out that we can obtain a suitable exceptional sequence which recovers the semiorthogonal decomposition after pushing down simply by mutating the existing sequence appropriately. 

Firstly, starting with the exceptional sequence 
\begin{equation}\label{eq:originalseq}
    \sod{\sh{E}_0, \sh{E}_1, \sh{F}_0, \sh{F}_1, \sh{A}_0, \sh{A}_1}
\end{equation}
we find that the objects $\sh{E}_1, \sh{F}_0$ are orthogonal in \Cref{subsec:morphisminres} below, so that 
\[
\sod{\sh{E}_0, \sh{F}_0, \sh{E}_1, \sh{F}_1, \sh{A}_0, \sh{A}_1}
\]
is also an exceptional sequence.
Then we perform a series of left and right mutations, using the distinguished triangles described in \Cref{sec:mutation}.
\begin{thm}\label{thm:resolution}
    Let
    \begin{align*}
        & \sh{A}_i^{\prime\prime} := \mathbb{L}_{\sh{E}_1}\mathbb{L}_{\sh{F}_1}\sh{A}_i \\
        & \sh{E}_0^{\prime\prime} := \mathbb{R}_{\sh{A}^{\prime\prime}_1}\mathbb{R}_{\sh{A}^{\prime\prime}_0}\sh{E}_0 \\
        & \sh{F}_0^{\prime\prime} := \mathbb{R}_{\sh{A}^{\prime\prime}_1}\mathbb{R}_{\sh{A}^{\prime\prime}_0}\sh{F}_0.
    \end{align*}
    Then there is a semiorthogonal deccomposition of $\db{\widetilde{X}}$ given by 
    \[
    \db{\widetilde{X}} = \sod{\sh{A}_0^{\prime\prime}, \sh{A}_1^{\prime\prime}, \sh{E}_0^{\prime\prime}, \sh{F}_0^{\prime\prime}, \sh{E}_1, \sh{F}_1, \db{B}}
    \]
    such that:
    \begin{enumerate}
        \item $\sigma_{*}\sh{A}_0^{\prime\prime} \cong \Phi(\sh{L}(1))$ and $\sigma_{*}\sh{A}_1^{\prime\prime} \cong \Phi(\sh{L}(2))$;
        \item the admissible subcategory $\widetilde{\sh{D}}$ defined using the exceptional sequence
        \[
        \widetilde{\sh{D}} := \sod{\sh{E}_0^{\prime\prime}, \sh{F}_0^{\prime\prime}, \sh{E}_1, \sh{F}_1}
        \]
        induces an equivalence
        \[
        \widetilde{\sh{D}} / (\ker \sigma_{*} \cap \widetilde{\sh{D}}) \xrightarrow{\sim} \sh{D}
        \]
        where $\sh{D} = \gens{\Phi(\sh{P})} \subset \db{X}$ is the absorbing category.
    \end{enumerate}
\end{thm}

The mutations performed in the above theorem are as follows:
\begin{align*}
    & \sod{\sh{E}_0, \sh{F}_0, \sh{E}_1, \sh{F}_1, \sh{A}_0, \sh{A}_1} \\
    & \curly \sod{\sh{E}_0, \sh{F}_0, \sh{E}_1, \sh{A}_0^{\prime}, \sh{A}_1^{\prime}, \sh{F}_1} \\
    & \curly \sod{\sh{E}_0, \sh{F}_0, \sh{A}_0^{\prime\prime}, \sh{A}_1^{\prime\prime}, \sh{E}_1, \sh{F}_1} \\
    & \curly \sod{\sh{A}_0^{\prime\prime}, \sh{A}_1^{\prime\prime}, \sh{E}_0^{\prime\prime}, \sh{F}_0^{\prime\prime}, \sh{E}_1, \sh{F}_1}
\end{align*}
where we firstly mutate $\sod{\sh{A}_0, \sh{A}_1}$ to the left of $\sod{\sh{E}_1, \sh{F}_1}$, and then mutate $\sod{\sh{E}_0, \sh{F}_0}$ to the right of $\sod{\sh{A}_0^{\prime\prime}, \sh{A}_1^{\prime\prime}}$.
The theorem then states that pushforward under $\sigma_{*}$ recovers the semiorthogonal decomposition of $\db{X}$.
In particular, once the image of each component under $\sigma_{*}$ is identified, the final statement follows as a consequence of Bondal-Orlov localization, as described in \Cref{sec:boloc}.
Since $\widetilde{\sh{D}}$ is an admissible component of a smooth projective variety, we suggest to think of it as a categorical resolution for the absorbing category.
Note that there are various formalisations of the notion of a categorical resolution in the literature, see for example \citep{kuz08} and \citep{kuznetsov-lunts}; here we use the term less formally.

In the following \Cref{sec:morphisms,sec:resmutation}, we firstly obtain a complete description of the morphisms in the original exceptional sequence \eqref{eq:originalseq}, as well the morphisms between the exceptional objects and to the objects in $\ker \sigma_{*}$.
Using this, we are then able to perform the necessary mutations, and give a proof of \Cref{thm:resolution}.

\subsection{Morphisms in the exceptional sequence}\label{sec:morphisms}

\subsubsection{Morphisms to the null category}
We firstly compute the morphisms between each object in the collection and some key objects in the kernel.
The results are summarised in \Cref{tab:nullhoms1} and \Cref{tab:nullhoms2}.

\begin{table}[!ht]
    \centering
    \begin{tabular}{ | c | c | c | c | }
    \hline
     \diagbox{$A$}{$B$} & $\sheaf_C(-1)$ & $\sheaf_E(-1, -1)$ & $\sheaf_E(-1, 0)$ \\ \hline
     $\sh{A}_i$ & $\mathbb{C}$ & 0 & 0 \\  \hline
     $\sh{F}_0$ & $\mathbb{C}[-1]$ & $\mathbb{C}$ & $\mathbb{C}^{\oplus 2}$ \\  \hline
     $\sh{F}_1$ & $\mathbb{C}[-1]$ & 0 & $\mathbb{C}$ \\  \hline
     $\sh{E}_0$ & 0 & $\mathbb{C}[-1]$ & $\mathbb{C}^{\oplus 2}[-1]$ \\ \hline
     $\sh{E}_1$ & 0 & 0 & $\mathbb{C}[-1]$ \\ \hline
    \end{tabular}
    \vspace{1.5em}
    \caption{$\Extcomplex(A, B)$ where $B$ is an object in the null category}
    \label{tab:nullhoms1}
\end{table}

\begin{table}[!ht]
    \centering
    \begin{tabular}{ | r | c | c | c | c | c | }
    \hline
     \diagbox{$A$}{$B$} & $\sh{A}_i$ & $\sh{F}_0$ & $\sh{F}_1$ & $\sh{E}_0$ & $\sh{E}_1$ \\ \hline
     $\sheaf_{C}(-1)$     & $\mathbb{C}[-3]$ & $\mathbb{C}[-2]$ & $\mathbb{C}[-2]$ & $0$ & $0$ \\ \hline
     $\sheaf_{E}(-1, -1)$ & $0$ & $\mathbb{C}[-2]$ & $\mathbb{C}^{\oplus 2}[-2]$ & $\mathbb{C}[-1]$ & $\mathbb{C}^{\oplus 2}[-1]$\\ \hline
     $\sheaf_{E}(-1, 0)$  & $0$ & $0$ & $\mathbb{C}[-2]$ & $0$ & $\mathbb{C}[-1]$ \\ \hline
    \end{tabular}
    \caption{$\Extcomplex(A, B)$ where $A$ is an object in the null category}
    \label{tab:nullhoms2}
\end{table}

\begin{enumerate}[wide, labelwidth=!, labelindent=0pt]
\item \textbf{Morphisms to $\sheaf_C(-1)$}

The objects $\exobj_0$ and $\exobj_1$ are orthogonal to $\sheaf_{C}(-1)$, because they are supported on $E_W$ which is disjoint from $C$.
Both $\sh{F}_0$ and $\sh{F}_1$ are supported on $E \cup E_W$, and we see that $E$ and $C$ intersect transversely, which gives us that:
\begin{align*}
\Extcomplex(\sheaff_i, \sheaf_{C}(-1)) \cong \mathbb{C}[-1]    
\end{align*}

Consider the objects $\sh{A}_i$, and let $j : \widetilde{F} \to \widetilde{X}$ be the inclusion.
We have a splitting distinguished triangle for $j^{*}\sh{A}_i \cong j^{*}j_{*}\sheaf(-1, i - 1)$:
\[
\sheaf(0, i - 1)[1] \to j^{*}\sh{A}_1 \to \sheaf(-1, i - 1)
\]
which shows that
\begin{align*}
    \Extcomplex(\sh{A}_i, \sheaf_C(-1)) & \cong \Extcomplex_{\widetilde{F}}(j^{*}\sh{A}_i, \sheaf_C(-1)) \\
    & \cong \Extcomplex_{\widetilde{F}}(\sheaf(-1, i - 1) \oplus \sheaf(0, i - 1)[1], \sheaf_C(-1)) \\
    & \cong \Extcomplex_{C}(\sheaf(-1) \oplus \sheaf[1], \sheaf(-1)) \cong \mathbb{C}.
\end{align*}

The above calculations also let us compute morphisms from $\sheaf_C(-1)$, since it is a $3$-spherical object by \Cref{lemma:threespherical}, so that $\Ext^{i}(\sheaf_C(-1), A) \cong \Ext^{3 - i}(A, \sheaf_C(-1))$ for any object $A$.

\item \textbf{Morphisms between $\sh{A}_i$ and $\sheaf_E(-1, a)$}

The objects $\sh{A}_0, \sh{A}_1$ are completely orthogonal to these objects because \[(\varphi_{Y})_{*}\sheaf_{E}(-1, a) = 0\] for any $a$, and because 
\[\Ext^{i}(\sheaf_{E}(-1, a), -) \cong \Ext^{3 - i}(-, \sheaf_{E}(-1, a - 2)).\]
by \Cref{restrictcanonical}.

\item \textbf{Morphisms between $\sh{F}_i$ and $\sheaf_E(-1, a)$}

Next we compute $\Extcomplex(\sheaff_0, \sheaf_{E}(-1, -1))$.
The key in this computation is that $E$ is the strict transform of $E_Y \subset W$, so:
\begin{align*}
    \Extcomplex_{\widetilde{X}}(\sheaff_0, \sheaf_{E}(-1, -1))
        & \cong \Extcomplex_{W}(\sheaf_{E_Y}(-1, -1), \sheaf_{E_Y}(-1, -1)) \cong \mathbb{C}
\end{align*}
and similarly
\begin{align*}
    \Extcomplex_{\widetilde{X}}(\sheaff_0, \sheaf_{E}(-1, 0))
        & \cong \Extcomplex_{W}(\sheaf_{E_Y}(-1, -1), \sheaf_{E_Y}(-1, 0)) \cong \mathbb{C}^{\oplus 2}
\end{align*}
because $\gens{\sheaf_{E_Y}(-1, -1), \sheaf_{E_Y}(-1, 0)} \simeq \db{\pp}$.
We similarly find that
\begin{align*}
    & \Extcomplex_{\widetilde{X}}(\sheaff_1, \sheaf_{E}(-1, 0))
        \cong \mathbb{C} \\
    & \Extcomplex_{\widetilde{X}}(\sheaff_1, \sheaf_{E}(-1, -1))
        \cong 0
\end{align*}
To compute morphisms in the reverse direction, we use Serre duality as before.
\begin{align*}
    & \Extcomplex(\sheaf_{E}(-1, -1), \sheaff_0) \cong \mathbb{C}[-2] \\
    & \Extcomplex(\sheaf_{E}(-1, -1), \sheaff_1) \cong \mathbb{C}^{\oplus 2}[-2] \\
    & \Extcomplex(\sheaf_{E}(-1, 0), \sheaff_0) = 0 \\
    & \Extcomplex(\sheaf_{E}(-1, 0), \sheaff_1) \cong \mathbb{C}[-2]
\end{align*}

\item \textbf{Morphisms between $\sh{E}_i$ and $\sheaf_E(-1, a))$}

Instead of using a geometric method to compute morphisms from $\sh{E}_i$, we make the following observation.
\begin{lemma}\label{lemma:sesei}
    For $i \in \{0, 1\}$, there are short exact sequences:
    \begin{equation}\label{eq:kr1triangle}
        0 \to \sh{E}_i \to \sh{F}_i \to \sheaf_{E}(-1, i - 1) \to 0.
    \end{equation}
\end{lemma}
\begin{proof}
We can form the following commuting diagram where the rows and first two columns are exact sequences.
\[\begin{tikzcd}
	& 0 & 0 & 0 \\
	0 & {\mathcal{O}} & {\mathcal{O}(E_W)} & {\mathcal{O}_{E_W}(-1, 0)} & 0 \\
	0 & {\mathcal{O}} & {\mathcal{O}(E + E_W)} & {\varphi_W^{*}\mathcal{O}_{E_Y}(-1, 0)} & 0 \\
	0 & 0 & {\mathcal{O}_E(-1, 0)} & {\mathcal{O}_E(-1, 0)} & 0 \\
	& 0 & 0 & 0
	\arrow[from=2-1, to=2-2]
	\arrow[from=2-2, to=2-3]
	\arrow[from=2-3, to=2-4]
	\arrow[from=2-4, to=2-5]
	\arrow[from=1-4, to=2-4]
	\arrow[from=2-4, to=3-4]
	\arrow[from=3-4, to=4-4]
	\arrow[from=4-4, to=5-4]
	\arrow[from=4-4, to=4-5]
	\arrow[from=3-4, to=3-5]
	\arrow[from=3-3, to=3-4]
	\arrow[from=3-2, to=3-3]
	\arrow[from=3-1, to=3-2]
	\arrow[from=1-2, to=2-2]
	\arrow[from=2-2, to=3-2]
	\arrow[from=3-2, to=4-2]
	\arrow[from=4-2, to=5-2]
	\arrow[from=4-2, to=4-3]
	\arrow[from=4-3, to=4-4]
	\arrow[from=3-3, to=4-3]
	\arrow[from=2-3, to=3-3]
	\arrow[from=1-3, to=2-3]
	\arrow[from=4-1, to=4-2]
	\arrow[from=4-3, to=5-3]
\end{tikzcd}\]
Applying the nine-lemma then gives us the desired short exact sequence from the rightmost column, and running the argument with the appropriate line bundle twist gives us the sequence for $\sh{E}_0$.
\end{proof}

The remaining morphisms in the tables are an immediate consequence of these distinguished triangles and the following statements:
\begin{align*}
    & \Extcomplex(\sheaf_{E}(a, b), \sheaf_{E}(a, b)) \cong \mathbb{C} \oplus \mathbb{C}[-2] \\
    & \Extcomplex(\sheaf_{E}(a, b), \sheaf_{E}(a, b - 1)) = 0
\end{align*}
which follow using the usual distinguished triangle \eqref{eq:counittriangle} for $\sheaf_{E}(a, b)\restrict{E}$:
\[
\sheaf(a + 2, b - 1)[1] \to \sheaf_{E}(a, b)\restrict{E} \to \sheaf(a, b)
\]
where we have used that $\sh{N}_{E / \widetilde{X}} \cong \sheaf(-2, 1)$, and invoked \cref{eq:linebundlepushforward} in order to compute cohomology of line bundles on the ruled surface $E$.

\end{enumerate}

\subsubsection{Morphisms between objects in the exceptional sequence}\label{subsec:morphisminres}

We summarise the results in \Cref{tab:sequencehoms}.

\begin{table}[!ht]
    \centering
    \begin{tabular}{ | c | c | c | c | c | c | c | }
    \hline
     \diagbox{$A$}{$B$} & $\sh{E}_0$ & $\sh{E}_1$ & $\sh{F}_0$ & $\sh{F}_1$ & $\sh{A}_0$ & $\sh{A}_1$ \\ \hline
     $\sh{E}_0$ & $\mathbb{C}$ & $\mathbb{C}^{\oplus 2}$ & $\mathbb{C} \oplus \mathbb{C}[-1]$ & $\mathbb{C}^{\oplus 2} \oplus \mathbb{C}^{\oplus 2}[-1]$ & $\mathbb{C}[-1]$ & $\mathbb{C}[-1]$ \\ \hline
     $\sh{E}_1$ & $0$ & $\mathbb{C}$ & $0$ & $\mathbb{C} \oplus \mathbb{C}[-1]$ & $\mathbb{C}[-1]$ & $\mathbb{C}[-1]$ \\ \hline
     $\sh{F}_0$ & $0$ & $0$ & $\mathbb{C}$ & $\mathbb{C}^{\oplus 2}$ & $\mathbb{C}[-1]$ & $\mathbb{C}[-1]$ \\ \hline
     $\sh{F}_1$ & $0$ & $0$ & $0$ & $\mathbb{C}$ & $\mathbb{C}[-1]$ & $\mathbb{C}[-1]$ \\ \hline
     $\sh{A}_0$ & $0$ & $0$ & $0$ & $0$ & $\mathbb{C}$ & $\mathbb{C}^{\oplus 2}$ \\ \hline
     $\sh{A}_1$ & $0$ & $0$ & $0$ & $0$ & $0$ & $\mathbb{C}$ \\ \hline
    \end{tabular}
    \caption{$\Extcomplex(A, B)$ where $A, B$ are objects in the exceptional sequence.}
    \label{tab:sequencehoms}
\end{table}

\begin{enumerate}[wide, labelindent=0pt]
\item \textbf{Morphisms to $\sh{A}_i$}

To begin, we see that
\begin{align*}
    \Ext^i_{\widetilde{X}}(\sheaff_1, \sh{A}_1)
        & \cong \Ext^i_{W}(\sheaf_{E_Y}(-1, 0), \psi^{*}\sheaf_{E_B}(-1, 0)) & \\
        & \cong \Ext^{3 - i}_{W}(\psi^{*}\sheaf_{E_B}(-1, 0), \sheaf_{E_Y}(-2, -1)) & \\
        & \cong \Ext^{3 - i}_{W}(\sheaf_{E_B}(-1, 0), \psi_{*}\sheaf_{E_Y}(-2, -1)) & \\
        & \cong \Ext^{3 - i}_{Y}(\sheaf_{E_B}(-1, 0), \sheaf_{Z_W}(-2)[1]) \cong \mathbb{C} & \text{for $i = 1$}
\end{align*}
and that there are no morphisms in any other degree.
One sees, since the map is induced by a transversal intersection at a single point, that the same result holds for $\Extcomplex(\sheaff_j, \sh{A}_k)$ for any choice of $j, k \in \{0, 1\}$.
Since the objects $\sh{A}_i$ are orthogonal to $\sheaf_{E}(-1, a)$ for any $a$, we find that $\Extcomplex(\sh{E}_j, \sh{A}_k) \cong \Extcomplex(\sh{F}_j, \sh{A}_k)$ for any $j, k \in \{0, 1\}$.

\item \textbf{Morphisms to $\sh{F}_i$} 

For each $i \in \{0, 1\}$, using that $\Extcomplex(\sh{E}_i, \sheaf_E(-1, i - 1)) \cong \mathbb{C}[-1]$, we obtain that \[\Extcomplex(\sh{E}_i, \sh{F}_i) \cong \mathbb{C} \oplus \mathbb{C}[-1].\]
Since $\Extcomplex(\sh{E}_1, \sheaf_E(-1, -1)) = \Extcomplex(\sh{E}_1, \sh{E}_0) = 0$, we must also have that $\Extcomplex(\sh{E}_1, \sh{F}_1) = 0$.
Finally, using that $\Extcomplex(\sh{E}_0, \sheaf_E(-1, 0)) \cong \mathbb{C}^{\oplus 2}[-1]$, we find that \[\Extcomplex(\sh{E}_0, \sh{F}_1) \cong \mathbb{C}^{\oplus 2} \oplus \mathbb{C}^{\oplus 2}[-1].\]

\end{enumerate}

\subsection{Mutating the semiorthogonal decomposition}\label{sec:resmutation}

We are now in a position to perform the mutations described previously.
The first set is defined by the following exact triangles.
\begin{align}
    & \sh{A}_0 \to \sh{A}_0^{\prime} \to \sh{F}_1 \label{eq:mut1a0} \\
    & \sh{A}_1 \to \sh{A}_1^{\prime} \to \sh{F}_1 \label{eq:mut1a1}
\end{align}
We have that the objects $\sh{A}_i^{\prime}$ are sheaves supported on $\widetilde{F} \cup E_W \cup E$.
Now note that the orthogonality of $\sh{A}_0, \sh{A}_1$ with $\sheaf_{E}(-1, 0)$ implies that there is an isomorphism:
\[
\Extcomplex(\sh{F}_1, \sh{A}_j) \cong \Extcomplex(\sh{E}_1, \sh{A}_j) \cong \mathbb{C}[-1]
\]
which we obtain by applying $\Extcomplex(-, \sh{A}_j)$ to the triangle \eqref{eq:kr1triangle}.
This isomorphism is obtained by pre-composing with $\sh{E}_1 \to \sh{F}_1$, so we obtain a commuting triangle.
\[\begin{tikzcd}
	{\mathcal{E}_1} & {\mathcal{F}_1} \\
	& {\mathcal{A}_j[1]}
	\arrow[from=1-2, to=2-2]
	\arrow[from=1-1, to=1-2]
	\arrow[from=1-1, to=2-2]
\end{tikzcd}\]
Applying $\Extcomplex(\sh{E}_1, -)$ to the triangles \eqref{eq:mut1a0} and \eqref{eq:mut1a1} gives us exact sequences:
\begin{align*}
    0 & \to \Hom(\sh{E}_1, \sh{F}_1) \xrightarrow{\sim} \Ext^1(\sh{E}_1, \sh{A}_j) \\ & \xrightarrow{0} \Ext^{1}(\sh{E}_1, \sh{A}^{\prime}_j) \xrightarrow{\sim} \Ext^1(\sh{E}_1, \sh{F}_1) \to 0
\end{align*}
where the first non-zero map is an isomorphism due to the above commuting triangle.
Hence, we have $\Ext^1(\sh{E}_1, \sh{A}_j^{\prime}) \cong \Ext^1(\sh{E}_1, \sh{F}_1) \cong \mathbb{C}$.
The next set of mutations is therefore defined by the following exact triangles.
\begin{align}
    & \sh{A}^{\prime}_0 \to \sh{A}_0^{\prime\prime} \to \sh{E}_1 \label{eq:mut2a0} \\
    & \sh{A}^{\prime}_1 \to \sh{A}_1^{\prime\prime} \to \sh{E}_1 \label{eq:mut2a1}
\end{align}

We now have the job of tidying up the exceptional sequence:
\begin{align*}
& \sod{\sh{E}_0, \sh{F}_0, \sh{A}_0^{\prime\prime}, \sh{A}_1^{\prime\prime}, \sh{E}_1, \sh{F}_1} \\
    & \curly \sod{\sh{A}_0^{\prime\prime}, \sh{A}_1^{\prime\prime}, \sh{E}_0^{\prime\prime}, \sh{F}_0^{\prime\prime}, \sh{E}_1, \sh{F}_1}    
\end{align*}

\begin{lemma}
    We have that
    \[
    \Extcomplex(\sh{F}_0, \sh{A}_0^{\prime\prime}) \cong \Extcomplex(\sh{F}_0, \sh{A}_1^{\prime\prime}) \cong \mathbb{C}.
    \]
    and if $\sheaff_0^{\prime} := \operatorname{Cone}(\sh{F}_0 \to \mutt{\sh{A}_0})$, then we have that
    \[
    \Extcomplex(\mut{\sh{F}}_0, \sh{A}_1^{\prime\prime}) \cong \mathbb{C}.
    \]
\end{lemma}
\begin{proof}
    Applying $\Extcomplex(\sh{F}_0, -)$ to \eqref{eq:mut2a0} gives an isomorphism $\Extcomplex(\sh{F}_0, \sh{A}_0^{\prime}) \cong \Extcomplex(\sh{F}_0, \sh{A}_0^{\prime\prime})$, since $\Extcomplex(\sh{F}_0, \sh{E}_1) = 0$. Similarly applying $\Extcomplex(\sh{F}_0, -)$ to \eqref{eq:mut1a0} gives a long exact sequence of $\Ext$ groups: \begin{align*}
      0 \to & \Hom(\sh{F}_0, \sh{A}_0) \to \Hom(\sh{F}_0, \sh{A}^{\prime}_0) \to \Hom(\sh{F}_0, \sh{F}_1) \to \\
      & \Ext^1(\sh{F}_0, \sh{A}_0) \to \Ext^1(\sh{F}_0, \sh{A}^{\prime}_0) \to \Ext^1(\sh{F}_0, \sh{F}_1) \to 0
\end{align*}
which reads:
\begin{align*}
    0 \to & 0 \to \Hom(\sh{F}_0, \sh{A}^{\prime}_0) \to \mathbb{C}^{\oplus 2} \to \\
      & \mathbb{C} \to \Ext^1(\sh{F}_0, \sh{A}^{\prime}_0) \to 0
\end{align*}
Then, the map $\mathbb{C}^{\oplus 2} \cong \Hom(\sh{F}_0, \sh{F}_1) \to \Ext^1(\sh{F}_0, \sh{A}_0) \cong \mathbb{C}$ is non-zero, due to the fact that the extension of $\sh{F}_0$ by $\sh{A}_0$ is induced by the intersection of $E_W \cup E$ with $\widetilde{F}$, and there is a map $\sh{F}_0 \to \sh{F}_1$ which is non-vanishing in a neighbourhood of the intersection locus, which is a reducible curve.

The second part of the lemma amounts to showing that the map \[\mathbb{C}^{\oplus 2} \cong \Hom(\mutt{\sh{A}_0}, \mutt{\sh{A}_1}) \to \Hom(\sh{F}_0, \mutt{\sh{A}_1}) \cong \mathbb{C}\] induced by pre-composing with $\sheaff_0 \to \mutt{\sh{A}_0}$ is non-zero.
By mutating over $\sh{E}_1$ and then $\sh{F}_1$, we see that this is equivalent to showing that the map \[
\mathbb{C}^{\oplus 2} \cong \Hom(\sh{A}_0, \sh{A}_1) \to \Ext^1(\sh{F}_2, \sh{A}_1) \cong \mathbb{C}
\]
is non-zero, where $\sh{F}_2 = \varphi_{W}^{*}\sheaf_{E_Y}(-1, 1)$ and the morphism is given by shifting and then precomposing with $\sh{F}_2 \to \sh{A}_1[1]$.
The argument to show that this map is non-zero is analagous to the one in the previous paragraph.
\end{proof}

We use the above lemma to create the following distinguished triangles defining an object $\sh{F}_0^{\prime\prime}$:
\begin{align}
    & \sh{F}_0 \to \sh{A}_0^{\prime\prime} \to \sh{F}_0^{\prime} \label{mutf01} \\
    & \sh{F}^{\prime}_0 \to \sh{A}_1^{\prime\prime} \to \sh{F}_0^{\prime\prime} \label{mutf02}
\end{align}
This object has a more explicit description.

\begin{lemma}
    Let $\sh{K}$ be the object defined by the extension
    \begin{equation} \label{eq:ktri}
        0 \to \sheaf_C(-1) \to \sh{K} \to \sheaf_E(-1, 0) \to 0.
    \end{equation}
    Then $\sh{F}_0^{\prime\prime}$ is given by the extension
    \begin{equation} \label{eq:f0ppe1}
        0 \to \sh{K} \to \sh{F}_0^{\prime\prime} \to \sh{E}_1 \to 0.
    \end{equation}
\end{lemma}
\begin{proof}
Since $\sh{E}_1$ is orthogonal to $\sheaf_C(-1)$, we have $\Extcomplex(\sh{E}_1, \sh{K}) \cong \mathbb{C}[-1]$ giving us an extension:
\[
0 \to \sh{K} \to \sh{G} \to \sh{E}_1 \to 0.
\]

It is clear by applying $\Extcomplex(\sh{E}_1, -)$ to this short exact sequence that $\Extcomplex(\sh{E}_1, \sh{G}) = 0$.
Next, we apply $\Extcomplex(\sh{F}_1, -)$ to the triangle \eqref{eq:ktri} defining $\sh{K}$. 
This gives an exact sequence:
\begin{align*}
    0 & \to \Hom(\sh{F}_1, \sh{K}) \to \Hom(\sh{F}_1, \sheaf_E(-1, 0)) \cong \mathbb{C} \\
    & \to \mathbb{C} \cong \Ext^{1}(\sh{F}_1, \sheaf_C(-1)) \to \Ext^1(\sh{F}_1, \sh{K}) \to 0. 
\end{align*}
The claim is that the middle map is an isomorphism, so that $\Extcomplex(\sh{F}_1, \sh{K}) = 0$.
But this is clear geometrically, since the extensions $\sh{F}_1 \to \sheaf_C(-1)[1]$ and $\sheaf_E(-1, 0) \to \sheaf_C(-1)[1]$ are induced by the intersection of $C$ and $E$, and in a neighbourhood of the intersection point $\sh{F}_1 \to \sheaf_E(-1, 0)$ is precisely the identity.
Hence $\Extcomplex(\sh{F}_1, \sh{K}) = \Extcomplex(\sh{F}_1, \sh{E}_1) = 0$, so we get that $\Extcomplex(\sh{F}_1, \sh{G}) = 0$.

Since there are no morphisms from $\sheaf_C(-1)$ or $\sheaf_E(-1, 0)$ to $\sh{E}_0$, we see that 
\[\Extcomplex(\sh{G}, \sh{E}_0^{\prime\prime}) = 0.\]
Hence $\sh{G} \in \gens{\sh{F}_0^{\prime\prime}}$. 
It is clear by the mutations relating $\sh{F}_0^{\prime\prime}$ to $\sh{F}_0$ that $\Extcomplex(\sh{F}_0, \sh{F}_0^{\prime\prime}) \cong \mathbb{C}$.
Hence it suffices to show that $\Extcomplex(\sh{F}_0, \sh{G}) \cong \mathbb{C}$.
We know that $\Extcomplex(\sh{F}_0, \sheaf_E(-1, 0)) \cong \mathbb{C}^{\oplus 2}$, and by similarly reasoning as we did to show that the map $\Hom(\sh{F}_1, \sheaf_E(-1, 0)) \to \Ext^1(\sh{F}_1, \sheaf_C(-1))$ is non-zero, we also find that the map
\[
\mathbb{C}^{\oplus 2} \cong \Hom(\sh{F}_0, \sheaf_E(-1, 0)) \to \Ext^1(\sh{F}_0, \sheaf_C(-1)) \cong \mathbb{C}
\]
is non-zero.
Hence, $\Extcomplex(\sh{F}_0, \sh{K}) \cong \mathbb{C}$, and since $\Extcomplex(\sh{F}_0, \sh{E}_1) = 0$ it follows immediately that $\Extcomplex(\sh{F}_0, \sh{G}) \cong \mathbb{C}$.
\end{proof}

Next we mutate the object $\sh{E}_0$.
\begin{lemma}
    We have that
    \[
    \Extcomplex(\sh{E}_0, \sh{A}_0^{\prime\prime}) \cong \Extcomplex(\sh{E}_0, \sh{A}_1^{\prime\prime}) \cong \mathbb{C}.
    \]
    and if $\mut{\sh{E}_0} := \operatorname{Cone}(\sh{E}_0 \to \sh{A}_0^{\prime\prime})$ then
    \[
    \Extcomplex(\sh{E}_0^{\prime}, \sh{A}_1^{\prime\prime}) \cong \mathbb{C}.
    \]
\end{lemma}
\begin{proof}
    The first statement holds because of the exact triangle:
    \[
    \sh{E}_0 \to \sh{F}_0 \to \sheaf_{E}(-1, -1)
    \]
    and the fact that $\mutt{\sh{A}_0}, \mutt{\sh{A}_1}$ are orthogonal to $\sheaf_{E}(-1, -1)$.
    For the second statement, we use the nine-lemma for derived categories, \Cref{lemma:derivednine}, and find that there is a distinguished triangle:
    \[
    \sh{E}_0^{\prime} \to \sh{F}_0^{\prime} \to \sheaf_{E}(-1, -1)[1]
    \]
\end{proof}

We similarly define the object $\sh{E}_0^{\prime\prime}$:
\begin{align}
    & \sh{E}_0 \to \sh{A}_0^{\prime\prime} \to \sh{E}_0^{\prime} \label{mute01} \\
    & \sh{E}^{\prime}_0 \to \sh{A}_1^{\prime\prime} \to \sh{E}_0^{\prime\prime} \label{mute02}
\end{align}
The cone on the morphism $\sh{E}_0^{\prime\prime} \to \sh{F}_0^{\prime\prime}$ is seen to be $\sheaf_E(-1, -1)[2]$ by another application of the derived nine-lemma.
In particular $\sh{E}_0^{\prime\prime}$ is not a sheaf.

Recall that we know that $\sigma_{*}\sh{E}_1 \cong \Phi(\sh{P})$, and an easy corollary of \Cref{lemma:sesei} is that $\sigma_{*}\sh{F}_1 \cong \sigma_{*}\sh{E}_1$.
The next lemmas will allow us to describe what happens if we pushforward the objects $\sh{A}_i^{\prime\prime}$.

\begin{lemma}\label{lemma:descAi}
    For each $j \in \{0, 1\}$, the objects $\sh{A}_j$ and $\mutt{\sh{A}}_j$ are related by the exact triangle
    \[
    \sh{A}_j \to \mutt{\sh{A}}_j \to \sh{Y}
    \]
    where $\sh{Y}$ is given by the non-trivial extension
    \[
        \sheaff_1 \to \sh{Y} \to \sh{E}_1.
    \]
\end{lemma}
\begin{proof}
    Consider the exact triangles
    \begin{align*}
        & \sh{A}_0 \to \sh{A}_0^{\prime} \to \sh{F}_1 \\
        & \sh{A}^{\prime}_0 \to \sh{A}_0^{\prime\prime} \to \sh{E}_1 \\
        & \sh{A}_0 \to \sh{A}_0^{\prime\prime} \to \sh{Y}
    \end{align*}
    Here, we use that $\Extcomplex(\sh{A}_0, \sh{A}_0^{\prime\prime}) \cong \Extcomplex(\sh{A}_0, \sh{A}_0^{\prime}) \cong \Extcomplex(\sh{A}_0, \sh{A}_0) \cong \mathbb{C}$, where the first isomorphism is induced by post-composing with $\sh{A}_0^{\prime} \to \sh{A}_0^{\prime\prime}$ and the second isomorphism is induced by post-composing with $\sh{A}_0 \to \sh{A}_0^{\prime}$.
    Therefore, the map $\sh{A}_0 \to \sh{A}_0^{\prime\prime}$ is a composition of the maps $\sh{A}_0 \to \sh{A}_0^{\prime}$ and $\sh{A}_0^{\prime} \to \sh{A}_0^{\prime\prime}$, and we can use the octahedral axiom to obtain an exact triangle:
    \[
    \sh{F}_1 \to \sh{Y} \to \sh{E}_1.
    \]
    The triangle is non-split since the morphism $\sh{E}_1 \to \sh{F}_1[1]$ in the triangle is the composition $\sh{E}_1 \to \sh{A}_0^{\prime}[1] \to \sh{F}_1[1]$, and we know that composing with $\sh{A}_0^{\prime} \to \sh{F}_1$ induces an isomorphism $\Ext^{1}(\sh{E}_1, \sh{A}_1^{\prime}) \cong \Ext^{1}(\sh{E}_1, \sh{F}_1)$.
\end{proof}

It is geometrically easier to see what happens if we first pushforward to $\widehat{X}$; the next lemma shows that the kernel of this contraction is simply generated by $\sheaf_E(-1, -1), \sheaf_E(-1, 0)$.
We can see this as a direct consequence of the fact that for a surface node the kernel of the contraction is generated by the pushforward of the $\sheaf(-1)$ bundle on the exceptional curve.

\begin{lemma}
    The kernel of the contraction $\widetilde{X} \to \widehat{X}$ is generated by
    \[ \ker \varphi_{\widehat{X}, *} = \gens{\sheaf_E(-1, -1), \sheaf_E(-1, 0)}.\]
\end{lemma}
\begin{proof}
    Let $Z \subset \widehat{X}$ be the reduced singular locus.
    Since the objects $\sheaf(-1, -1), \sheaf(-1, 0)$ generate the kernel of the contraction $E \to Z$, we want to show that the kernel of $\widetilde{X} \to \widehat{X}$ is generated by the pushforwards of these objects.
    We use \citep[Theorem 5.2]{kuznetsov-shinder} to show this, which requires us to show that
    \[
    \varphi_{\widehat{X}, *}\sheaf_{\widetilde{X}}(-mE) \cong \sh{J}^{m}_{Z}
    \]
    where $\sh{J}_{Z}$ is the ideal sheaf of $Z$.
    We can work in a neighbourhood of $Z$, but locally near $Z$ we find that $\widehat{X}$ is just $V \times \pp$ for a nodal surface $V$ and $\widetilde{X}$ is locally around $E$ a pullback of the resolution $\widetilde{V}$ of $V$, i.e. we have a pullback square.
    \[\begin{tikzcd}
	{\widetilde{V}\times \mathbb{P}^1} & {\widetilde{V}} \\
	{V\times\mathbb{P}^1} & V
	\arrow["\pi", from=1-2, to=2-2]
	\arrow["f", from=1-1, to=1-2]
	\arrow["\pi"', from=1-1, to=2-1]
	\arrow["f"', from=2-1, to=2-2]
	\arrow["\lrcorner"{anchor=center, pos=0.125}, draw=none, from=1-1, to=2-2]
\end{tikzcd}\]
    Let $\pi : \widetilde{V} \to V$ be the contraction, $F$ the exceptional curve and $Z_V$ the nodal point.
    The desired result follows from pulling back the corresponding equivalence for $V$, which was shown in \citep[Corollary 5.6]{kuznetsov-shinder}:
    \[
    \pi_{*}\sheaf_{\widetilde{V}}(-mF) \cong \sh{J}_{Z_V}^{m}
    \]
    and using flat base change:
    \[
    f^{*}\pi_{*} \simeq \pi_{*}f^{*}
    \]
    as needed.
\end{proof}

Consider the extension on $\db{\curve}$:
\[
0 \to \sh{P} \to \sh{Q} \to \sh{P} \to 0
\]
so that $\sh{Q}$ is the skyscraper sheaf along the non-reduced branch $\pp_{u:v} \times \mathbb{C}[s]/s^2$.

\begin{lemma}\label{lemma:pushforwardY}
    We have that
    \begin{align*}
        & \sigma_{*}\sh{Y} \cong \Phi(\sh{Q})
    \end{align*}
\end{lemma}
\begin{proof}
    Firstly, we see that $\sh{Y}$ is orthogonal to the objects $\sheaf_{E}(-1, a)$ for all $a \in \mathbb{Z}$.
    Since \[\Extcomplex(\sh{E}_1, \sheaf_{E}(-1, -1)) = \Extcomplex(\sh{F}_1, \sheaf_E(-1, -1)) = 0\] we immediately have $\Extcomplex(\sh{Y}, \sheaf_E(-1, -1)) = 0$.
       We also know that \[\Ext^{1}(\sh{E}_1, \sh{F}_1) \cong \Ext^{1}(\sh{E}_1, \sheaf_E(-1, 0)),\] i.e. the map $\sh{E}_1 \to \sheaf_E(-1, 0)[1]$ factors through $\sh{E}_1 \to \sh{F}_1[1]$.
       Then, applying the functor $\Extcomplex(-, \sheaf_E(-1, 0))$ to the triangle $\sh{F}_1 \to \sh{Y} \to \sh{E}_1$ gives us an exact sequence:
       \begin{align*}
           0 & \to \Hom(\sh{Y}, \sheaf_E(-1, 0)) \to \Hom(\sh{F}_1, \sheaf_E(-1, 0)) \\
           & \xrightarrow{\sim} \Ext^1(\sh{E}_1, \sheaf_E(-1, 0) \to \Ext^{1}(\sh{Y}, \sheaf_E(-1, 0) \to 0
       \end{align*}
       which gives us that $\Extcomplex(\sh{Y}, \sheaf_E(-1, 0)) = 0$.
       It follows then that $\Extcomplex(\sh{Y}, \sheaf_E(-1, a)) = 0$ for all integers $a$, and hence by Serre duality that $\Extcomplex(\sheaf_E(-1, a), \sh{Y}) = 0$ also.
       
       Now, we pushfoward to $\widehat{X}$.
       Here, we have a distinguished triangle:
       \[
       \varphi_{\widehat{X}, *}\sh{F}_1 \to \varphi_{\widehat{X}, *}\sh{Y} \to \varphi_{\widehat{X}, *}\sh{E}_1
       \]
       One sees that $\varphi_{\widehat{X}, *}\sh{F}_1 \cong \varphi_{\widehat{X}, *}\sh{E}_1$ is the $\sheaf(-1, 0)$ bundle on the strict transform of $G_{\red}$.
       The extension is non-split, since we showed that $\sh{Y}$ is orthogonal to $\ker \varphi_{\widehat{X}, *}$, so for example we know that $\varphi_{\widehat{X}, *}\sh{Y}$ has no endomorphisms in non-zero degrees. 
       Hence, $\varphi_{\widehat{X}, *}\sh{Y}$ is the $\sheaf(-1, 0)$ bundle on the strict transform of $G$.
       Now we can see that $\sigma_{*}\sh{Y}$ is the $\sheaf(-1, 0)$ bundle on $G$, or equivalently, $\Phi(\sh{Q})$.
\end{proof}

\begin{proof}[Proof of \Cref{thm:resolution}]
To prove that 
\begin{align*}
        & \sigma_{*}\sh{A}_0^{\prime\prime} \cong \Phi(\sh{L}(1)) \\
        & \sigma_{*}\sh{A}_1^{\prime\prime} \cong \Phi(\sh{L}(2))
\end{align*}
recall the defining triangle for $\sh{A}_1^{\prime}$ \eqref{eq:mut1a1}:
    \[
    \sh{A}_1 \to \sh{A}_1^{\prime} \to \sh{F}_1.
    \]
    In fact, we note that the transverse intersection of $E \cup E_W$ with $\widetilde{F}$ along $C_E \cup C_W \subset \widetilde{F}$ induces an extension between $\sh{F}_1$ and $\sh{A}_1$.
    By uniqueness, we find that $\sh{A}_1^{\prime}$ is the pushforward of a line bundle on $E \cup E_W \cup \widetilde{F}$ which is $\sheaf(-1, 0) \otimes \sheaf(C_E + C_W)$ on $\widetilde{F}$, $\sheaf(-1, 0)$ on $E$ and $\sheaf$ on $E_W$.
    One verifies that this is a valid line bundle on the union since all restrictions to intersections agree.
There is a line bundle on $\pp_{s:t} \cup \pp_{u:v} \subset C$ (i.e. the union of the reduced branches) which restricts to $\sheaf(1)$ on $\pp_{s:t}$ and $\sheaf$ on $\pp_{u:v}$.
We denote by $\sh{M}$ the pushforward of this into $\curve$.
Then one sees that $\sigma_{*}\sh{A}_1^{\prime} \cong \Phi(\sh{M})$.
Let $\sheaf_{\curvebase}$ be the skyscraper sheaf along the reduced branch.
Then $\sigma_{*}\sh{A}_1 = \Phi(\sheaf_{\curvebase})$.

Next, there is a unique extension in $\db{\curve}$ between $\sh{P}$ and $\sh{M}$:
\begin{gather}
    \sh{M} \to \sh{L}(2) \to \sh{P}. \label{eq:pushAtri2}
\end{gather}
We can consider the triangle obtained by pushing forward:
\begin{gather*}
    \sh{A}^{\prime}_1 \to \sh{A}_1^{\prime\prime} \to \sh{E}_1
\end{gather*}
If this pushforward is non-split, it corresponds to the triangle obtained by embedding \eqref{eq:pushAtri2} using $\Phi$.
To show that it does not split, we pushforward the triangle
\[
\sh{A}_1 \to \sh{A}_1^{\prime\prime} \to \sh{Y}
\]
which must correspond to an extension in $\db{\curve}$ between $\sh{Q}$ and $\sheaf_{\curvebase}$ by \Cref{lemma:pushforwardY}.
Then one computes that there is indeed an extension induced by intersection:
\[
\sheaf_{\curvebase} \to \sh{L}(2) \to \sh{Q}
\]
so that $\sigma_{*}\sh{A}_1^{\prime\prime}$ can only be $\Phi(\sh{L}(2))$.
A similar argument shows $\sigma_{*}\sh{A}_0^{\prime\prime} \cong \Phi(\sh{L}(1))$.

Now let
\[
\widetilde{\sh{D}} = \sod{\sh{E}_0^{\prime\prime}, \sh{F}_0^{\prime\prime}, \sh{E}_1, \sh{F}_1}.
\]
The statement
\[
\widetilde{\sh{D}}/(\ker \sigma_{*} \cap \widetilde{\sh{D}}) \xrightarrow{\sim} \sh{D}
\]
follows from the fact that the cones on the morphisms $\sh{F}_0^{\prime\prime} \to \sh{E}_1$, $\sh{E}_0^{\prime\prime} \to \sh{F}_0^{\prime\prime}$ and $\sh{E}_1 \to \sh{F}_1$ are in the kernel of $\sigma_{*}$, so the image of $\widetilde{\sh{D}}$ under $\sigma_{*}$ is contained within our absorbing category $\sh{D}$.
Combined with the fact that $\sh{A}_0^{\prime\prime}, \sh{A}_1^{\prime\prime}$ pushfoward to $\Phi(\sh{L}(1))$ and $\Phi(\sh{L}(2)$ respectively, we find that the image of $\widetilde{\sh{D}}$ under $\sigma_{*}$ is in fact the entire absorbing category $\sh{D}$, and the desired equivalence follows from Bondal-Orlov localization (see \Cref{sec:boloc}).
\end{proof}

\subsection{Structure of the categorical resolution}

Having found the subcategory $\widetilde{\sh{D}} \subset \db{\widetilde{X}}$, we now endeavour to describe the morphisms between the objects and how they behave under composition.
Since the subcategory is found by mutating an exceptional sequence, the resulting objects also form an exceptional sequence
\[
    \widetilde{\sh{D}} = \sod{\sh{E}_0^{\prime\prime}, \sh{F}_0^{\prime\prime}, \sh{E}_1, \sh{F}_1}.
\]
We recall that the objects $\sh{E}_1, \sh{F}_1$ are sheaves supported on the exceptional loci $E_W, \widetilde{E}_Y$ as defined in \cref{def:exceptionals}, while $\sh{E}_0^{\prime\prime}$ and $\sh{F}_0^{\prime\prime}$ are defined by the triangles induced by mutation given in \cref{mute01,mute02} and \cref{mutf01,mutf02} respectively.

\subsubsection{Morphisms in the categorical resolution}

\begin{thm}\label{thm:resolutionmorphisms}
    The morphisms in $\widetilde{\sh{D}}$ are described as follows:
    \begin{enumerate}
        \item $\Extcomplex(\sh{E}_0^{\prime\prime}, \sh{F}_0^{\prime\prime}) \cong \Extcomplex(\sh{E}_1, \sh{F}_1) \cong \mathbb{C} \oplus \mathbb{C}[-1]$
        \item $\Extcomplex(\sh{F}_0^{\prime\prime}, \sh{E}_1) \cong \mathbb{C} \oplus \mathbb{C}[-1]$
        \item $\Extcomplex(\sh{F}_0^{\prime\prime}, \sh{F}_1) \cong \Extcomplex(\sh{E}_0^{\prime\prime}, \sh{E}_1) \cong \mathbb{C} \oplus \mathbb{C}[-1] \oplus \mathbb{C}^{\oplus 2}[-2]$
        \item $\Extcomplex(\sh{E}_0^{\prime\prime}, \sh{F}_1) \cong \mathbb{C} \oplus \mathbb{C}[-1] \oplus \mathbb{C}^{\oplus 2}[-2] \oplus \mathbb{C}^{\oplus 2}[-3]$
    \end{enumerate}
\end{thm}

The first claim of the theorem has already been proven in \Cref{subsec:morphisminres}, and we prove the rest in a series of lemmas.

\begin{lemma} \label{lemma:morphismsg1}
    We have that $\Extcomplex(\sh{F}_0^{\prime\prime}, \sh{E}_1) \cong \mathbb{C} \oplus \mathbb{C}[-1]$ and $\Extcomplex(\sh{F}_0^{\prime\prime}, \sh{F}_1) \cong \mathbb{C} \oplus \mathbb{C}[-1] \oplus \mathbb{C}^{\oplus 2}[-2]$.
\end{lemma}
\begin{proof}
    Apply $\Extcomplex(-, \sh{E}_1)$ to \eqref{eq:f0ppe1}.
Since $\Extcomplex(\sheaf_C(-1), \sh{E}_1) = 0$, we find that 
\[\Extcomplex(\sh{K}, \sh{E}_1) \cong \Extcomplex(\sheaf_E(-1, 0), \sh{E}_1) \cong \mathbb{C}[-1]\] and this gives us that $\Extcomplex(\sh{F}_0^{\prime\prime}, \sh{E}_1) \cong \mathbb{C} \oplus \mathbb{C}[-1]$.

    Applying $\Extcomplex(-, \sh{F}_1)$ to the triangle \eqref{eq:ktri} defining $\sh{K}$, we find that $\Extcomplex(\sh{K}, \sh{F}_1) \cong \mathbb{C}^{\oplus 2}[-2]$ because we have
    \[
    \Extcomplex(\sheaf_C(-1), \sh{F}_1) \cong \Extcomplex(\sheaf_E(-1, 0), \sh{F}_1) \cong \mathbb{C}[-2].
    \]
    Then, applying $\Extcomplex(-, \sh{F}_1)$ to the triangle $\eqref{eq:f0ppe1}$ gives us that 
    \[\mathbb{C} \cong \Ext^{i}(\sh{E}_1, \sh{F}_1) \cong \Ext^i(\sh{F}_0^{\prime\prime}, \sh{F}_1)\]
    for $i = 0, 1$, and
    \[\Ext^2(\sh{F}_0^{\prime\prime}, \sh{F}_1) \cong \Ext^2(\sh{K}, \sh{F}_1) \cong \mathbb{C}^{\oplus 2}.\]
\end{proof}

\begin{lemma} \label{lemma:g0e1}
    We have that
    \begin{align*}
        & \Extcomplex(\sh{E}_0^{\prime\prime}, \sh{E}_1) \cong \mathbb{C} \oplus \mathbb{C}[-1] \oplus \mathbb{C}^{\oplus 2}[-1] \\
    & \Extcomplex(\sh{E}_0^{\prime\prime}, \sh{F}_1) \cong \mathbb{C} \oplus \mathbb{C}[-1] \oplus \mathbb{C}^{\oplus 2}[-2] \oplus \mathbb{C}^{\oplus 2}[-3].
    \end{align*}
\end{lemma}
\begin{proof}
    Applying $\Extcomplex(-, \sh{E}_1)$ to the exact triangle
    \[
    \sh{E}^{\prime\prime}_0 \to \sh{F}_0^{\prime\prime} \to \sheaf_{E}(-1, -1)[2]
    \]
    and using the fact that $\Extcomplex(\sheaf_{E}(-1, -1)[2], \sh{E}_1) \cong \mathbb{C}^{\oplus 2}[-3]$ we get the first statement.
    The second statement is proven by instead applying $\Extcomplex(-, \sh{F}_1)$, but now using that \[\Extcomplex(\sheaf_{E}(-1, -1)[2], \sh{F}_1) \cong \mathbb{C}
    ^{\oplus 2}[-4].\]
\end{proof}

\subsubsection{Algebra structure of the categorical resolution}
We now investigate how the compositions behave in the categorical resolution.
Notice that for a $\cpinf$ object $P$ we have
\[\Extcomplex(P, P) \cong \mathbb{C} \oplus \mathbb{C}[-1] \oplus \mathbb{C}^{\oplus 2}[-2] \oplus \mathbb{C}^{\oplus 2}[-3] \oplus \cdots\]
and that we can recover the vector spaces of morphisms between objects in our categorical resolution $\widetilde{\sh{D}}$ as truncations of this algebra.
\begin{align*}
    & \Extcomplex(\sh{E}_0^{\prime\prime}, \sh{F}_0^{\prime\prime}) \cong \Extcomplex(\sh{E}_1, \sh{F}_1)  \cong \Extcomplex(\sh{F}_0^{\prime\prime}, \sh{E}_1) \cong \Ext^{\leq 1}(P, P) && \cong \mathbb{C} \oplus \mathbb{C}[-1]\\
    & \Extcomplex(\sh{E}_0^{\prime\prime}, \sh{E}_1)  \cong \Extcomplex(\sh{F}_0^{\prime\prime}, \sh{F}_1) \cong \Ext^{\leq 2}(P, P) && \cong \mathbb{C} \oplus \mathbb{C}[-1] \oplus \mathbb{C}^{\oplus 2}[-2]\\
    & \Extcomplex(\sh{E}_0^{\prime\prime}, \sh{F}_1) \cong \Ext^{\leq 3}(P, P) && \cong \mathbb{C} \oplus \mathbb{C}[-1] \oplus \mathbb{C}^{\oplus 2}[-2] \oplus \mathbb{C}^{\oplus 2}[-3]
\end{align*}
Subsequently, for any two objects $A, B$ in the exceptional sequence, we think of an element of $\Ext^{i}(A, B)$ as lifting some element of $\Ext^{i}(P, P)$. 
Since $\Extcomplex(P, P)$ is a polynomial ring with no relations, this suggests the following.

\begin{prop}\label{prop:algebrastructure}
    All compositions in $\widetilde{\sh{D}}$ are non-zero.
\end{prop}

Before proving the proposition, we illustrate some consequences.
Let $t_0, t_1, t_2, s_0, s_1$ be non-zero degree zero maps and $\epsilon_0, \epsilon_1, \epsilon_2, \delta_0, \delta_1$ be non-zero degree one maps as described by the following quiver.
    \begin{equation}\label{quiver:morphisms}
    \begin{tikzcd}[row sep=large,column sep=large]
    	{\mathcal{E}_0^{\prime\prime}} && {\mathcal{F}_0^{\prime\prime}} \\
    	\\
    	{\mathcal{E}_1} && {\mathcal{F}_1}
    	\arrow["{t_0}"{description}, curve={height=-6pt}, from=1-1, to=1-3]
    	\arrow["{\epsilon_0}"{description}, curve={height=6pt}, dashed, from=1-1, to=1-3]
    	\arrow["{\delta_0}"{description}, curve={height=6pt}, dashed, from=1-1, to=3-1]
    	\arrow["{s_0}"{description}, curve={height=-6pt}, from=1-1, to=3-1]
    	\arrow["{t_1}"{description}, curve={height=6pt}, from=1-3, to=3-1]
    	\arrow["{\epsilon_1}"{description}, curve={height=-6pt}, dashed, from=1-3, to=3-1]
    	\arrow["{t_2}"{description}, curve={height=-6pt}, from=3-1, to=3-3]
    	\arrow["{\epsilon_2}"{description}, curve={height=6pt}, dashed, from=3-1, to=3-3]
    	\arrow["{s_1}"{description}, curve={height=6pt}, from=1-3, to=3-3]
    	\arrow["{\delta_1}"{description}, curve={height=-6pt}, dashed, from=1-3, to=3-3]
    \end{tikzcd}
    \end{equation}
Letting $ab$ denote the composition $b \circ a$, a corollary of \Cref{prop:algebrastructure} is that up to rescaling, we can obtain the following relations
\begin{align*}
    & t_0t_1 = s_0 && t_1t_2 = s_1\\    
    & t_0\epsilon_1 = \epsilon_0t_1 = \delta_0 && t_1\epsilon_2 = \epsilon_1t_2 = \delta_1.
\end{align*}
Defining $\rho_0 := \epsilon_0 \epsilon_1$ and $\rho_1 := \epsilon_1 \epsilon_2$, we also see that 
\begin{align*}
    & \rho_0t_2 = t_0 \rho_1 && \rho_0\epsilon_2 = \epsilon_0 \rho_1.
\end{align*}
We can then extend to find bases $\{\rho_0, \theta_0\}$ of $\Ext^{2}(\sh{E}_0^{\prime\prime}, \sh{E}_1)$ and $\{\rho_1, \theta_1\}$ of $\Ext^{2}(\sh{F}_0^{\prime\prime}, \sh{F}_1)$ respectively.\footnote{Note that we can choose $\theta_0, \theta_1$ such that either of the relations $t_0\theta_1 = \theta_0t_2$ or $\epsilon_0 \theta_1 = \theta_0 \epsilon_2$ hold.
However, it is not clear whether we can choose $\theta_0, \theta_1$ such that \textit{both} commutativity relations simultaneously hold.}
We can then visualise the category $\widetilde{\mathcal{D}}$ through the following quiver, where we omit compositions in non-zero degrees, the solid arrows are maps of degree $0$, the dashed arrows are maps in degree $1$ and the dotted arrows are maps in degree $2$.
    \[\begin{tikzcd}[row sep=large,column sep=large]
	{\mathcal{E}_0^{\prime\prime}} && {\mathcal{F}_0^{\prime\prime}} \\
	\\
	{\mathcal{E}_1} && {\mathcal{F}_1}
	\arrow["{t_0}"{description}, curve={height=-6pt}, from=1-1, to=1-3]
	\arrow["{\epsilon_0}"{description}, curve={height=6pt}, dashed, from=1-1, to=1-3]
	\arrow["{s_0}"{description}, curve={height=-6pt}, from=1-1, to=3-1]
	\arrow["{t_1}"{description}, curve={height=6pt}, from=1-3, to=3-1]
	\arrow["{\epsilon_1}"{description}, curve={height=-6pt}, dashed, from=1-3, to=3-1]
	\arrow["{\theta_1}"{description}, curve={height=-6pt}, dotted, from=1-3, to=3-3]
	\arrow["{t_2}"{description}, curve={height=-6pt}, from=3-1, to=3-3]
	\arrow["{\epsilon_2}"{description}, curve={height=6pt}, dashed, from=3-1, to=3-3]
	\arrow["{s_1}"{description}, curve={height=6pt}, from=1-3, to=3-3]
	\arrow["{\theta_0}"{description}, curve={height=6pt}, dotted, from=1-1, to=3-1]
\end{tikzcd}\]
In particular, we highlight how our resolving category is essentially built out of the graded Kronecker quivers (see \Cref{def:gradedkronecker}) that resolve $\mathbb{P}^{\infty}$ objects.

The proof of \Cref{prop:algebrastructure} comprises the remainder of the section through a series of lemmas.

\begin{enumerate}[wide, labelindent=0pt]
\item \textbf{Compositions with source $\sh{F}_0^{\prime\prime}$}

\begin{lemma}
    The compositions
    \begin{align*}
    t_1t_2 & : \sh{F}_0^{\prime\prime} \to \sh{E}_1 \to \sh{F}_1 \\ 
    t_1\epsilon_2 & : \sh{F}_0^{\prime\prime} \to \sh{E}_1 \to \sh{F}_1[1]
    \end{align*}
    are non-zero.
\end{lemma}
\begin{proof}
The first statement follows from the proof of \Cref{lemma:morphismsg1}, which shows that the degree zero map $\sh{F}_0^{\prime\prime} \to \sh{F}_1$ factors as $\sh{F}_0^{\prime\prime} \to \sh{E}_1 \to \sh{F}_1$. 
The degree $1$ map $\sh{F}_0^{\prime\prime} \to \sh{F}_1[1]$ similarly factors as $\sh{F}_0^{\prime\prime} \to \sh{E}_1 \to \sh{F}_1[1]$.
\end{proof}

\begin{lemma}
    The composition $\epsilon_1t_2 : \sh{F}_0^{\prime\prime} \to \sh{F}_1$ is non-zero.
\end{lemma}
\begin{proof}    
We deduce that $\Extcomplex(\sh{F}_0^{\prime\prime}, \sheaf_E(-1, 0)) \cong \mathbb{C}^{\oplus 2}[-2]$, so we can immediately see that the composition $\sh{F}_0^{\prime\prime} \to \sh{E}_1[1] \to \sh{F}_1[1]$ is non-zero.
\end{proof}

\begin{lemma}
    The composition $\epsilon_1\epsilon_2 : \sh{F}_0^{\prime\prime} \to \sh{E}_1[1] \to \sh{F}_1[2]$ is non-zero.
\end{lemma}
\begin{proof}
    Consider the following diagram.
\begin{equation} \label{eq:ninediagram2}
\begin{tikzcd}
	{\mathcal{K}} & {\mathcal{O}_E(-1, 0)} & {\mathcal{O}_C(-1)[1]} \\
	{\mathcal{F}_0^{\prime\prime}} & {\mathcal{E}_1[1]} & {\mathcal{Y}^{\prime}[1]} \\
	{\mathcal{E}_1} & {\mathcal{F}_1[1]} & {\mathcal{Y}^{\prime\prime}}
	\arrow[from=1-1, to=2-1]
	\arrow[from=1-1, to=1-2]
	\arrow[from=1-2, to=1-3]
	\arrow[from=1-2, to=2-2]
	\arrow[from=2-2, to=3-2]
	\arrow[from=2-1, to=2-2]
	\arrow[from=2-1, to=3-1]
	\arrow[dashed, from=3-1, to=3-2]
	\arrow[dashed, from=3-2, to=3-3]
	\arrow[from=2-2, to=2-3]
	\arrow[dashed, from=1-3, to=2-3]
	\arrow[dashed, from=2-3, to=3-3]
\end{tikzcd}
\end{equation}
There is a unique degree $1$ morphism (up to scaling) $\sh{K} \to \sh{E}_1[1]$, and one can quickly show that this factors through both $\sh{K} \to \sh{F}_0^{\prime\prime}$ and $\sh{K} \to \sheaf_E(-1, 0)$, by applying $\Extcomplex(-, \sh{E}_1)$ to the triangles defining $\sh{K}$ and $\sh{F}_0^{\prime\prime}$ respectively.
Hence, by the derived nine-lemma we can complete along the dashed arrows to obtain exact triangles.
The bottom left square must then commute, so it follows that the map $\sh{E}_1 \to \sh{F}_1[1]$ is non-zero and hence we find that $\sh{Y}^{\prime\prime} \cong \sh{Y}[1]$.
Then, we have an exact triangle
\[
\sheaf_C(-1) \to \sh{Y}^{\prime} \to \sh{Y}.
\]
We see that $\Hom(\sh{Y}^{\prime}, \sh{F}_1) \cong \Hom(\sh{Y}, \sh{F}_1) \cong \mathbb{C}$ and $\Ext^2(\sh{Y}^{\prime}, \sh{F}_1) \cong \Ext^2(\sheaf_C(-1), \sh{F}_1) \cong \mathbb{C}$, and there are no morphisms in any other degrees.
Then, applying $\Extcomplex(-, \sh{F}_1)$ to the triangle
\[
\sh{E}_1 \to \sh{Y}^{\prime} \to \sh{F}_0^{\prime\prime}
\]
we see that since $\Ext^1(\sh{Y}^{\prime}, \sh{F}_1) = 0$ that the map $\Ext^1(\sh{E}_1, \sh{F}_1) \to \Ext^2(\sh{F}_0^{\prime\prime}, \sh{F}_1)$ is an injection and hence the composition $\sh{F}_0^{\prime\prime} \to \sh{E}_1[1] \to \sh{F}_1[2]$ is non-zero.
\end{proof}

\item \textbf{Compositions with source $\sh{E}_0^{\prime\prime}$ and target $\sh{E}_1$}

\begin{lemma}
    The compositions 
    \begin{align*}
        t_0t_1 & : \sh{E}_0^{\prime\prime} \to \sh{F}_0^{\prime\prime} \to \sh{E}_1  \\
        \epsilon_0t_1 & : \sh{E}_0^{\prime\prime} \to \sh{F}_0^{\prime\prime}[1] \to \sh{E}_1[1] \\
        t_0\epsilon_1 & : \sh{E}_0^{\prime\prime} \to \sh{F}_0^{\prime\prime} \to \sh{E}_1[1]
    \end{align*}
    are non-zero.
\end{lemma}
\begin{proof}
    Consider the exact triangle
\[
\sh{K} \to \sh{F}_0^{\prime\prime} \to \sh{E}_1
\]
and apply the $\Extcomplex(\sh{E}_0^{\prime\prime}, -)$ functor. 
Using the fact that $\sh{E}_0^{\prime\prime}$ is orthogonal to $\sheaf_C(-1)$ and the defining exact triangle for $\sh{K}$
\[
\sheaf_C(-1) \to \sh{K} \to \sheaf_E(-1, 0)
\]
we find that \[\Extcomplex(\sh{E}_0^{\prime\prime}, \sh{K}) \cong \Extcomplex(\sh{E}_0^{\prime\prime}, \sheaf_E(-1, 0)) \cong \mathbb{C}^{\oplus 2}[-3].\]
It follows that there is an isomorphism
\[
\Hom(\sh{E}_0^{\prime\prime}, \sh{F}_0^{\prime\prime}) \cong \Hom(\sh{E}_0^{\prime\prime}, \sh{E}_1) \cong \mathbb{C}
\]
so that the composition $\sh{E}_0^{\prime\prime} \to \sh{F}_0^{\prime\prime} \to \sh{E}_1$ is non-zero, and similarly:
\[
\Ext^1(\sh{E}_0^{\prime\prime}, \sh{F}_0^{\prime\prime}) \cong \Ext^1(\sh{E}_0^{\prime\prime}, \sh{E}_1) \cong \mathbb{C}
\]
so the composition $\sh{E}_0^{\prime\prime} \to \sh{F}_0^{\prime\prime}[1] \to \sh{E}_1[1]$ is non-zero.
We can also consider the exact triangle
\[
\sh{E}_0^{\prime\prime} \to \sh{F}_0^{\prime\prime} \to \sheaf_E(-1, -1)[2]
\]
and apply the $\Extcomplex(-, \sh{E}_1)$ functor.
Using that $\Extcomplex(\sheaf_E(-1, -1)[2], \sh{E}_1) \cong \mathbb{C}^{\oplus 2}[-3]$, we find an isomorphism:
\[
\Ext^1(\sh{E}_0^{\prime\prime}, \sh{E}_1) \cong \Ext^1(\sh{F}_0^{\prime\prime}, \sh{E}_1) \cong \mathbb{C}
\]
which implies that the composition $\sh{E}_0^{\prime\prime} \to \sh{F}_0^{\prime\prime} \to \sh{E}_1[1]$ is non-zero.
\end{proof}

\begin{lemma}
    The composition $\epsilon_0\epsilon_1 : \sh{E}_0^{\prime\prime} \to \sh{F}_0^{\prime\prime}[1] \to \sh{E}_1[2]$ is non-zero.
\end{lemma}
\begin{proof}
Suppose for a contradiction that the composition is zero. 
Applying $\Extcomplex(-, \sh{F}_0^{\prime\prime})$ to the triangle \[\sh{E}_0^{\prime\prime} \to \sh{F}_0^{\prime\prime} \to \sheaf_E(-1, -1)[2]\] gives us an exact sequence:
\begin{align*}
0 = & \Ext^1(\sh{F}_0^{\prime\prime}, \sh{F}_0^{\prime\prime}) \to \Ext^1(\sh{E}_0^{\prime\prime}, \sh{F}_0^{\prime\prime}) \to \\ & \Ext^2(\sheaf_E(-1, -1)[2], \sh{F}_0^{\prime\prime}) \to \Ext^2(\sh{F}_0^{\prime\prime}, \sh{F}_0^{\prime\prime}) = 0
\end{align*}
so that we have an isomorphism $\mathbb{C} \cong \Ext^1(\sh{E}_0^{\prime\prime}, \sh{F}_0^{\prime\prime}) \to \Ext^2(\sheaf_E(-1, -1)[2], \sh{F}_0^{\prime\prime}) \cong \mathbb{C}$.
Hence the left triangle commutes in the following diagram.
\[\begin{tikzcd}
	{\mathcal{O}_E(-1, -1)[1]} & {\mathcal{E}_0^{\prime\prime}} \\
	{\mathcal{F}_0^{\prime\prime}[1]} & {\mathcal{E}_1[2]}
	\arrow[from=1-1, to=1-2]
	\arrow[from=1-1, to=2-1]
	\arrow[from=1-2, to=2-1]
	\arrow[from=2-1, to=2-2]
	\arrow["0", from=1-2, to=2-2]
\end{tikzcd}\]
Since the right triangle commutes by assumption, it suffices to show that the composition 
\begin{equation}\label{eq:compg1e1}
\sheaf_E(-1, -1) \to \sheaff_0^{\prime\prime} \to \sh{E}_1[1]    
\end{equation}
is non-zero, since then we obtain a contradiction.
For this, we consider the non-split extension
\[
\sh{E}_1 \to \sh{Y}^{\prime} \to \sheaff_0^{\prime\prime}
\]
and we only need to show that $\Hom(\sheaf_E(-1, -1), \sh{Y}^{\prime}) = 0$, since then the map 
\[\Hom(\sheaf_E(-1, -1), \sheaff_0^{\prime\prime}) \to \Ext^1(\sheaf_E(-1, -1), \sh{E}_1)\]
is an injection, and hence the composition \eqref{eq:compg1e1} is non-zero.
Recall the diagram \eqref{eq:ninediagram2} gives us an exact triangle
\[
\sheaf_C(-1) \to \sh{Y}^{\prime} \to \sh{Y}.
\]
There are no morphisms from $\sheaf_{E}(-1, -1)$ to $\sh{E}_1$ and $\sh{F}_1$ in degree $0$.
Hence \[\Hom(\sheaf_E(-1, -1), \sh{Y}) = 0\] and since $\Hom(\sheaf_E(-1, -1), \sheaf_C(-1)) = 0$, we have our desired result.
\end{proof}

\item \textbf{Compositions with source $\sh{E}_0^{\prime\prime}$ and target $\sh{F}_1$}
\begin{lemma}
    For all $i \in \{0, 1, 2\}$, precomposing with $t_0 : \sh{E}_0^{\prime\prime} \to \sh{F}_0^{\prime\prime}$ induces isomorphisms
    \[
    \Ext^i(\sh{E}_0^{\prime\prime}, \sh{F}_1) \cong \Ext^i(\sheaff_0^{\prime\prime}, \sh{F}_1)
    \]
    and postcomposing with $t_2 : \sh{E}_1 \to \sh{F}_1$ induces isomorphisms
    \[
    \Ext^i(\sh{E}_0^{\prime\prime}, \sh{F}_1) \cong \Ext^i(\sh{E}_0^{\prime\prime}, \sh{E}_1).
    \]
\end{lemma}
\begin{proof}
Firstly, applying $\Extcomplex(-, \sh{F}_1)$ to the triangle
\[
\sh{E}_0^{\prime\prime} \to \sh{F}_0^{\prime\prime} \to \sheaf_E(-1, -1)[2]
\]
induces isomorphisms $\Ext^i(\sheaff_0^{\prime\prime}, \sh{F}_1) \cong \Ext^i(\sh{E}_0^{\prime\prime}, \sh{F}_1)$ by pre-composing with $\sh{E}_0^{\prime\prime} \to \sh{F}_0^{\prime\prime}$, for $i \in \{0, 1, 2\}$.
Similarly, using that $\Extcomplex(\sh{E}_0^{\prime\prime}, \sheaf_E(-1, 0)) \cong \mathbb{C}^{\oplus 2}[-3]$ and applying $\Extcomplex(\sh{E}_0^{\prime\prime}, -)$ to the triangle \eqref{eq:kr1triangle}
\[
\sh{E}_1 \to \sh{F}_1 \to \sheaf_E(-1, 0)
\]
gives us isomorphisms
\[
\Ext^{i}(\sh{E}_0^{\prime\prime}, \sh{E}_1) \cong \Ext^i(\sh{E}_0^{\prime\prime}, \sh{F}_1)
\]
for $i \in \{0, 1, 2\}$ by composing with $\sh{E}_1 \to \sh{F}_1$.    
\end{proof}

\begin{lemma}
    Precomposing with $\epsilon_0 : \sh{E}_0^{\prime\prime} \to \sh{F}_0^{\prime\prime}[1]$ induces an isomorphism
    \[
    \Ext^3(\sh{E}_0^{\prime\prime}, \sh{F}_1) \cong \Ext^{2}(\sheaff_0^{\prime\prime}, \sh{F}_1)
    \]
    and postcomposing with $\epsilon_2 : \sh{E}_1 \to \sh{F}_1[1]$ induces an isomorphism
    \[
    \Ext^3(\sh{E}_0^{\prime\prime}, \sh{F}_1) \cong \Ext^{2}(\sh{E}_0^{\prime\prime}, \sh{E}_1).
    \]
\end{lemma}
\begin{proof}
By applying $\Extcomplex(\sh{F}_0^{\prime\prime}, -)$ to 
\begin{equation} \label{eq:extendef}
    \sh{F}_1 \to \sh{Y} \to \sh{E}_1
\end{equation}
and using the results above, we find $\Extcomplex(\sh{F}_0^{\prime\prime}, \sh{Y}) \cong \mathbb{C} \oplus \mathbb{C}[-2]$.
Then, we apply the functor $\Extcomplex(\sheaf_E(-1, -1)[2], -)$ to the same triangle and see that $\Extcomplex(\sheaf_E(-1, -1)[2], \sh{Y})$ can only be non-zero in degrees $0 \leq i \leq 3$.
Hence, applying $\Extcomplex(-, \sh{Y})$ to the triangle
\[
\sh{E}_0^{\prime\prime} \to \sh{F}_0^{\prime\prime} \to \sheaf_E(-1, -1)[2]
\]
we find that $\Ext^3(\sh{E}_0^{\prime\prime}, \sh{Y}) = 0$.
By applying $\Extcomplex(\sh{E}_0^{\prime\prime}, -)$ to \eqref{eq:extendef}
we obtain an isomorphism
\[
\mathbb{C}^{\oplus 2} \cong \Ext^2(\sh{E}_0^{\prime\prime}, \sh{E}_1) \xrightarrow{\sim} \Ext^3(\sh{E}_0^{\prime\prime}, \sh{F}_1) \cong \mathbb{C}^{\oplus 2}.
\]
This shows the first statement.

We use a similar argument for the second statement.
Consider the non-split extension
\[
\sh{F}_0^{\prime\prime} \to \sh{Z}_0 \to \sh{E}_0^{\prime\prime}
\]
and apply $\Extcomplex(-, \sh{F}_1)$.
We see that if $\Ext^3(\sh{Z}_0, \sh{F}_1) = 0$ then there is an isomorphism
\[
\mathbb{C}^{\oplus 2} \cong \Ext^2(\sh{F}_0^{\prime\prime}, \sh{F}_1) \to \Ext^3(\sh{E}_0^{\prime\prime}, \sh{F}_1)  \cong \mathbb{C}^{\oplus 2}.
\]
We use the following diagram where the top left square commutes, and so we complete the dashed arrows to exact triangles using the derived nine-lemma.
\[\begin{tikzcd}
	{\mathcal{E}_0^{\prime\prime}} & {\mathcal{F}_0^{\prime\prime}[1]} & {\mathcal{Z}_0[1]} \\
	{\mathcal{F}_0^{\prime\prime}} & {\mathcal{E}_1[1]} & {\mathcal{Y}^{\prime}[1]} \\
	{\mathcal{O}_E(-1, -1)[2]} & {\mathcal{K}[2]} & {\mathcal{Z}_1[1]}
	\arrow[from=1-1, to=1-2]
	\arrow[from=1-2, to=1-3]
	\arrow[from=1-1, to=2-1]
	\arrow[from=2-1, to=3-1]
	\arrow[from=1-2, to=2-2]
	\arrow[from=2-1, to=2-2]
	\arrow[dashed, from=3-1, to=3-2]
	\arrow[from=2-2, to=3-2]
	\arrow[from=2-2, to=2-3]
	\arrow[dashed, from=3-2, to=3-3]
	\arrow[dashed, from=1-3, to=2-3]
	\arrow[dashed, from=2-3, to=3-3]
\end{tikzcd}\]
One quickly sees that $\Extcomplex(\sheaf_E(-1, -1)[2], \sh{F}_1)$ and $\Extcomplex(\sh{K}[2], \sh{F}_1)$ are only non-zero in degree $3$, and hence $\Ext^i(\sh{Z}_1, \sh{F}_1) = 0$ for $i \geq 4$.
From the previous calculations, we have that $\Extcomplex(\sh{Y}^{\prime}, \sh{F}_1) \cong \mathbb{C} \oplus \mathbb{C}[-2]$.
Hence, we find that $\Ext^3(\sh{Z}_0, \sh{F}_1) = 0$, by using the triangle 
\[
\sh{Z}_0 \to \sh{Y}^{\prime} \to \sh{Z}_1.
\]
\end{proof}

\end{enumerate}

\section{Deformation absorption through \texorpdfstring{$\cpinf$}{compound P-infinity} objects} \label{section:smoothing}

In this section we return to the singular curve $\curve$ defined in \cref{eq:curve}, and briefly discuss an interesting example of how $\cpinf$ objects behave when smoothing the curve.
Our goal is to investigate whether these objects provide deformation absorptions, as described in \Cref{subsec:defabs}.
For this section, we recall the semiorthogonal decomposition \eqref{eq:sodforC}:
\[
\db{\curve} = \sod{\sh{P}(-1), \db{\pp_{s:t}} }.
\]
We will denote $P := \sh{P}(-1)$ for this section.
For our purposes, we will define a smoothing $\mathfrak{X}$ of $\curve$ to be a flat family of curves $\phi : \mathfrak{X} \to (B, o)$, where $(B, o)$ is a smooth pointed curve, the central fibre is $\phi^{-1}(o) \cong \curve$ and all other fibres are smooth.
Note that we don't require the total space $\mathfrak{X}$ to be smooth.

Consider the surface 
    \[
        \mathfrak{X} = \bl{x, s^2}\pp_{x:y} \times \pp_{s:t}.
    \]
with the natural map $\mathfrak{X} \to \pp_{x:y}$.
We can embed into a projective bundle:
\[
\mathfrak{X} \cong \{ xv = s^2u \} \subset \pbundle_{\pp \times \pp}(\sheaf(0, 2)_u \oplus \sheaf(1, 0)_v)
\]
so that the space $\mathfrak{X}$ is nodal, and the fibre over $\{x = 0\}$ is given by \[ \{s^2u = 0\} \subset \pbundle_{\pp_{s:t}}(\sheaf(2)_u \oplus \sheaf_v)\]
which is isomorphic to the curve $\curve$.
Note that the generic fibre is a copy of $\pp_{s:t}$, so that $\mathfrak{X}$ is a smoothing of $\curve$.

\begin{prop}\label{prop:defabs}
    The object $P \in \db{\curve}$ provides a deformation absorption with respect to the smoothing $\mathfrak{X}$.
\end{prop}
\begin{proof}
    The blow up formula shows that there is a semiorthogonal decomposition
    \[
    \db{\mathfrak{X}} = \sod{\db{\mathbb{C}[s]/s^2}, \db{\pp \times \pp}}.
    \]
    If $i : \curve \xhookrightarrow{} \mathfrak{X}$ is the inclusion, then $i_{*}P$ is precisely the $\pinf{1}$ object generating the admissible subcategory $\db{\mathbb{C}[s]/s^2}$ of $\db{\mathfrak{X}}$.
\end{proof}

It is not clear whether the $\cpinf$ object provides a \textit{universal} deformation absorption, i.e. whether the pushforward of $P$ generates an admissible category for all smoothings of $\curve$.
We note that in our example, the pushforward of the $\cpinf$ object $P$ gives a $\pinf{1}$ object on $\mathfrak{X}$, which absorbs the singularities of the total space of the smoothing.
One observes that the total space of a smoothing of $\curve$ cannot itself be smooth, for purely homological reasons.
Indeed, consider the object $i_{*}P$.
The triviality of the normal bundle of $\curve$ in $\mathfrak{X}$ gives us the distinguished triangle using \cref{eq:counittriangle}:
\[
i^{*}i_{*}P \to P \to P[2].
\]
In the case that $P \to P[2]$ is the zero map, it is clear that $i_{*}P$ is not in $\dperf{\mathfrak{X}}$ since $i^{*}i_{*}P$ is not in $\dperf{\curve}$.
Otherwise, applying $\Extcomplex(-, P)$ to the triangle shows that
\[
\Ext^{j}(i_{*}P, i_{*}P) \cong \Ext^{j}(i^{*}i_{*}P, P) \cong \mathbb{C}
\]
for all $j \geq 0$, which implies $i_{*}P$ is not perfect.
Hence, for another smoothing, we may similarly require that $i_{*}P$ absorbs the singularities of the total space.

\DeclareFieldFormat[misc]{title}{\mkbibquote{#1}}
\printbibliography

\end{document}